\documentclass{amsart}

\usepackage{graphicx}
\usepackage{amssymb}
\usepackage{epstopdf}
\usepackage{float}
\usepackage{comment}
\usepackage{soul}

\usepackage{dsfont}
\usepackage{subcaption}
\usepackage{caption}
\usepackage{hyperref}
\usepackage[nameinlink,noabbrev]{cleveref}

\theoremstyle{theorem}
\newtheorem{theorem}[subsection]{Theorem}
\newtheorem{lemma}[subsection]{Lemma}

\theoremstyle{definition}
\newtheorem{definition}[subsection]{Definition}

\crefname{equation}{}{}
\Crefname{equation}{Equation}{Equations}
\crefname{figure}{Figure}{Figures}

\newcommand{\A}{\mathds{A}}
\newcommand{\C}{\mathds{C}}
\newcommand{\G}{\mathds{G}}
\newcommand{\I}{\mathds{I}}
\newcommand{\PP}{\mathds{P}}
\newcommand{\R}{\mathds{R}}
\newcommand{\U}{\mathds{U}}

\newcommand{\cH}{\mathcal{H}}

\newcommand{\cX}{\mathcal{X}}

\newcommand{\bp}{\textbf{p}}
\newcommand{\bu}{\textbf{u}}
\newcommand{\bv}{\textbf{v}}
\newcommand{\bw}{\textbf{w}}
\newcommand{\bx}{\textbf{x}}

\newcommand{\bK}{\textbf{K}}

\newcommand{\bT}{\textbf{T}}
\newcommand{\bW}{\textbf{W}}
\newcommand{\bX}{\textbf{X}}
\newcommand{\bY}{\textbf{Y}}

\newcommand{\ve}{\varepsilon}
\newcommand{\Id}{\textrm{Id}}
\newcommand{\imag}[1]{\textrm{Im} \left( #1 \right) }
\newcommand{\real}[1]{\textrm{Re} \left( #1 \right) }
\newcommand{\sige}[1]{\sigma_c\left( #1 \right) }

\newcommand{\diff}[2]{\frac{\mathrm{d} #1}{\mathrm{d} #2}}

\title{(In)Stability of Travelling Waves \\ in a Model of Haptotaxis \footnote{\today}}

\author[Harley et al.]{K.E. Harley, P. van Heijster, R. Marangell, \\ G.J. Pettet, T. V. Roberts and M. Wechselberger}

\begin{document}

\begin{abstract}
We examine the spectral stability of travelling waves of the haptotaxis model studied in \cite{Harley_2014}. In the process we apply Li\'enard coordinates to the linearised stability problem and use a Riccati-transform/Grassmanian spectral shooting method \`a la \cite{harley2015numerical,ledoux2009computing,ledoux2010grassmannian} in order to numerically compute the Evans function and point spectrum of a linearised operator associated with a travelling wave.  We numerically show the instability of non-monotone waves (type IV) and the stability of the monotone ones (types I-III) to perturbations in an appropriately weighted space.
\end{abstract}

\maketitle

\section[Introduction]{Introduction}\label{sec:set}

We study the system of partial differential equations (PDEs) introduced in \cite{PER} to describe haptotactic cell invasion in a model for melanoma. Haptotaxis, similar to chemotaxis, describes the preferred motion of cells towards, or away from, the gradient of a chemical concentration. This chemical is bound to a surface for haptotaxis, while it is suspended in a fluid for chemotaxis \cite{Harley_2014}.
The original proposed model in \cite{PER} considered three densities: the extracellular matrix (ECM) concentration, the invasive tumour cell population, and the density of protease. However, as the protease reaction was assumed to happen on a (super-)fast time scale \cite{PER}, a quasi-steady state approximation was used to reduce to a simplified model considering only the densities of the ECM and the tumour. Written in the nondimensionalised form of \cite{Harley_2014} that emphasises its advection-reaction-diffusion structure, the model is given by
\begin{equation}\label{eq:ardpde}
\begin{split}
\begin{pmatrix} u \\ w \end{pmatrix}_t & = \ve \begin{pmatrix} u \\ w \end{pmatrix}_{xx} + \begin{pmatrix} 0 \\ -wu_x \end{pmatrix}_x + \begin{pmatrix} -u^2w \\ w(1-w) \end{pmatrix}, 
\end{split}
\end{equation}
where $u$ and $w$ represent nondimensionalised concentrations of the ECM and the invasive tumour cell population respectively, and with $x \in \mathbb{R}, t \in \mathbb{R}^+$ and $\ve  \geq 0$ a small parameter\footnote{Note that the original model in \cite{PER} ignored diffusion ($\ve=0$) as it was assumed that diffusion only played a minimal role.}. 

\begin{figure}
	\makebox[\linewidth][c]{%
			\includegraphics[width=0.75\textwidth]{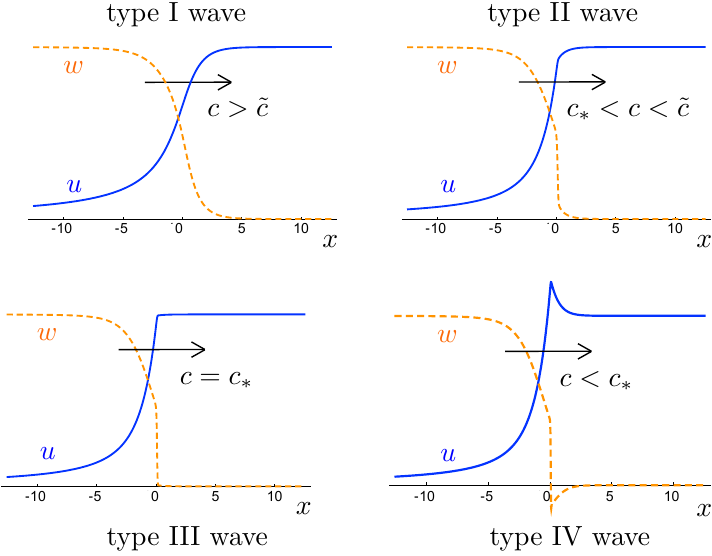}}
				\captionsetup{width=.95\linewidth}
	\caption{The four different types of travelling wave solutions supported by \eqref{eq:ardpde}.
	}\vspace{-.75cm}
	\label{fig:waves223}
\end{figure}

In \cite{Harley_2014} it was shown in a rigorous fashion that \eqref{eq:ardpde} supports four types of travelling wave solutions. 
The classification of the travelling wave solutions was based on distinguishing, qualitative features of the waves in the singular limit $\ve \to 0$.  {\em{Type I}} waves are smooth with a monotone wave profile, {\em{type II}} waves are shock-fronted in $w$ (in the singular limit $\ve\to 0$) with a monotone wave profile, {\em{type III}} waves are shock-fronted in $w$ with a monotone wave profile whose $w$-component has semi-compact support,
and {\em{type IV}} waves are shock-fronted in $w$ with a non-monotone wave profile (i.e. $w$ is negative for certain parts of the profile). \cref{fig:waves223} provides an example of the four types of waves found.

To arrive at this result \cite{Harley_2014} followed the work of \cite{Wechselberger_Pettet_10} and the model was analysed in its singular limit $\ve \to 0$ using canard theory and Li\'enard coordinates. Smooth travelling wave solutions (type I) were explicitly found for speeds larger than some critical speed $\tilde{c}$. Similarly, shock-fronted travelling wave solutions (type II-IV) were found for speeds smaller than this critical speed $\tilde{c}$. In particular, type II waves exist for speeds in between the so-called minimal wave speed $c_*$ \cite{PER} and the critical wave speed $\tilde{c}$, while type III waves travel with the minimal wave speed $c_*$ and type IV waves travel slower than the minimal wave speed $c_*$.
These travelling wave solutions were shown to persist for a small $\ve$ through the application of Geometric Singular Perturbation Theory (GSPT). 
These results extended/formalised the earlier results of \cite{HOS,PER}. 

The connection between the observed wave speed and the asymptotic behaviour of its initial condition was also investigated numerically in \cite{Harley_2014}. 
However, the (spectral) stability of these four types of travelling waves has not been determined before. Biologically, type IV waves are expected to be unstable simply because they contain regions with negative cell population. Furthermore, in \cite{marchant2006biphasic} it is argued that Type III waves are physically the most realistic as they have (i) sharp interfaces and (ii) zero tumour concentration in ahead of the interface. We numerically find that these waves correspond to stable waves with the smallest positive wave speed and that waves with smaller speeds (type IV waves) are unstable. Mathematically, the type III waves decay much faster at $+\infty$ than the type II or IV waves. This means that their derivative still decays in the appropriate exponentially weighted space. Hence, the temporal eigenvalue $\lambda =0$, associated with translation invariance, persists. This eigenvalue is a (locally) smooth function of the wave speed parameter $c$ and moves into the right-half plane as the wave speed is further decreased (as we numerically show). 

\subsection{Main result: spectral stability of type I-III waves and instability of type IV waves}

We numerically establish the stability of waves of type I-III and the instability of waves of type IV in appropriately exponentially weighted spaces via determination of the roots of an {\em Evans function}. Originally used in the determination of stability of nerve-axon impulses, Evans functions have received a boost in the last 30 years by linking stability of a travelling wave to geometric ideas \cite{agj90,allen2002numerical, beckmal13, carter2016stability, DGK2, DGK1,garzum98,HZ06, ledoux2010grassmannian, ledoux2009computing,Sandstede_02}. Computing the Evans function can be numerically delicate, and there are several geometrically inspired techniques to resolve this in the literature, \cite{agj90,allen2002numerical,bridges2002stability,de2016spectra,DGK1,gardner1991stability,grudzien2016geometric,ledoux2009computing,ledoux2010grassmannian}, to name a few. For a nice exposition of some of these as well as further references, see \cite{ledoux2010grassmannian}.

For our stability results, we will work on the Grassmannian as in \cite{ledoux2009computing,ledoux2010grassmannian}. 
The linearity of the spectral problem means it will induce a nonlinear flow on the Grassmannian \cite{brockett1981multivariable,lafortune1996superposition,ledoux2009computing,ledoux2010grassmannian,martin1978applications,schneider1973global,Shayman86}. Rather than keeping track of solutions themselves, since subspaces of solutions are preserved, we instead track them on the Grassmannian under the induced flow \cite{lafortune1996superposition,ledoux2009computing,ledoux2010grassmannian,martin1978applications,schneider1973global,Shayman86}. 
The flow induced by a linear system on the Grassmannian is called the {\emph{generalised (or extended) Riccati flow}} \cite{Shayman86}. It is a nonlinear, but lower order, flow on the manifold. 
The original definition of the Evans function can now be interpreted in terms of this Riccati flow on the Grassmannian, equivalently either through projection from the Steifel manifold \cite{ledoux2010grassmannian} onto a chart of the Grassmannian, or (as we do in this manuscript) via a meromorphic function which has been called the {\emph{Riccati-Evans function}} \cite{harley2015numerical}. Importantly, the solutions to the matrix Riccati equation seem to be numerically well behaved on the (charts of the) Grassmannian and we no longer have exponential growth of solutions \cite{ledoux2009computing,ledoux2010grassmannian}, though at the expense of some solutions becoming singular \cite{levin1959matrix}.   

Our evolution of the boundary data follows the Evans function calculation techniques developed in \cite{ledoux2009computing,ledoux2010grassmannian}, however, we have managed (in this case at least) to avoid the singularities which are typically present in solutions to the Riccati equation \cite{levin1959matrix,Shayman86}. 

Previous uses of the Riccati equation to generate an Evans function include \cite{de2016spectra,harley2015numerical,ledoux2009computing,ledoux2010grassmannian}. In \cite{ledoux2010grassmannian}, the Riccati-Evans function approach was used to confirm stability of Boussinesq solitary waves, autocatalytic travelling waves and the Ekman boundary layer. In \cite{ledoux2009computing}, the authors focussed on the stability of wrinkled fronts in a cubic autocatalysis reaction-diffusion system with two spatial independent variables. In \cite{de2016spectra}, the singular nature of the problem was exploited and used to generate a matrix Riccati equation and subsequent flow on the Grassmannian in order to study the stability of periodic pulse wavetrains. In \cite{harley2015numerical} the Riccati-Evans function approach was used to study the stability of travelling waves in two lower-dimensional models: the Fisher/Kolmogorov-Petrovsky-Piscounov equation and a Keller-Segel model of bacterial chemotaxis. In \cite{ledoux2009computing,ledoux2010grassmannian}, a chart changing mechanism was described to avoid singularities of the Riccati equation on the fly, and the method was linked to the so-called `continuous orthogonalisation' method \cite{HZ06,ledoux2009computing}, while in \cite{harley2015numerical} it was observed that by carefully picking a single standard chart, singularities could be avoided. 

The current manuscript shows another way to avoid singularities in the spectral parameter regime of interest. In particular, we do not work in the standard charts of the Grassmannian as in \cite{harley2015numerical}, but rather a judiciously chosen one.

This manuscript is organised as follows, 
in \cref{sec:back} we briefly discuss the key results of \cite{Harley_2014} needed for the stability analysis. 
In \cref{sec:spec} we describe the linearised problem and compute the essential and absolute spectrum of type I-IV waves. In \cref{sec:point} we expound on the Riccati-Evans function approach for computing the point spectrum and in \cref{sec:application} apply it to the haptotaxis model \eqref{eq:ardpde} to show the spectral instability of the type IV waves, as well as numerical evidence of spectral stability of waves of type I, II and III. In \cref{sec:fut} we briefly discuss related future research directions, both for the haptotaxis model \eqref{eq:ardpde} and the Riccati-Evans function.

\section{Setup: existence of travelling waves} \label{sec:back}
We reproduce the key results of \cite{Harley_2014} related to the existence of the four different types of travelling wave solutions (in a slightly modified form from \cite{Harley_2014}). 
Passing to a moving coordinate frame, we set $z = x-ct$ where $c>0$ is our wave speed parameter. We get the travelling wave form of the equation:
\begin{equation}\label{eq:twpde}
\begin{split}
\begin{pmatrix} u \\ w \end{pmatrix}_t & = \ve \begin{pmatrix} u \\ w \end{pmatrix}_{zz} + \begin{pmatrix} cu  \\ cw - wu_z \end{pmatrix}_z + \begin{pmatrix} -u^2w \\w(1-w)\end{pmatrix}. 
\end{split}
\end{equation}
A {\em travelling wave} will be a steady state solution to \cref{eq:twpde}, connecting two distinct background states of \cref{eq:ardpde}. The background states of \cref{eq:ardpde} are $(u,w) = (0,1)$ and $(u,w) = (u_\infty,0)$, for $u_\infty \in \R$ (i.e.\ we have a line of fixed points in \cref{eq:twode}). Thus, a travelling wave is a solution to the nonlinear ordinary differential equation (ODE) and in what follows we set $' := \frac{d}{dz}$ for notational convenience:
\begin{equation}\label{eq:twode}
\begin{split}
0 & = \ve \begin{pmatrix} u \\ w \end{pmatrix}'' + \begin{pmatrix} cu  \\ cw - wu' \end{pmatrix}' + \begin{pmatrix} -u^2w \\ w(1-w) \end{pmatrix}
\end{split}
\end{equation}
satisfying the boundary conditions 
\begin{equation}\label{eq:bcs}
\lim_{z \to -\infty}{u(z)} = 0, \quad \lim_{z \to +\infty}{u(z)} = u_{\infty}, \quad \lim_{z \to -\infty}{w(z)} = 1, \quad \lim_{z \to +\infty}{w(z)} = 0. 
\end{equation}
The second condition in \cref{eq:bcs} implies that the righthand boundary condition on $u$, denoted $u_\infty$ is free. In what follows we assume $u_\infty >0$. Introducing the variables (Li\'enard coordinates):
\begin{equation}\label{eq:lienardnonlin}
\begin{aligned}
v &:= u' \\
y &:= \ve w' - vw + cw
\end{aligned}
\end{equation}
allows us to re-write \cref{eq:twode} as a system of ODE with two fast ($v$ and $w$) and two slow ($u$ and $y$) variables: 
\begin{equation} \label{eq:nl-slow}
\begin{aligned}
u ' & =  v, \\
y' & =   -w(1-w), \\ 
\ve v' & =  -cv + u^2w, \\
\ve w' & =  y + vw - cw\,.
\end{aligned}
\end{equation}
We will refer to \cref{eq:nl-slow} as the {\em (nonlinear) slow system}, and the variable $z$ as the {\em slow travelling wave coordinate}. To investigate the problem in the fast timescale, we introduce the {\em fast travelling wave coordinate} $\zeta = z/\ve$ and derive the corresponding four dimensional {\em (nonlinear) fast system} with $\ve \neq 0$ and with the convention that $\dot{} := \frac{d}{d \zeta}$
\begin{equation}\label{eq:nl-fast}
\begin{aligned}
\dot{u} & =   \ve v, \\
\dot{y} & =  -\ve w (1-w), \\ 
\dot{v} & =  -cv + u^2 w, \\
\dot{w} & =  y + vw - cw\,.
\end{aligned}
\end{equation}
As in \cite{Harley_2014} we now set $\ve =0$ and pick out our solutions from the resulting systems. As $\ve \to 0$ the nonlinear fast system becomes the so-called {\em layer problem}
\begin{equation} \label{eq:nonlin-layer}
\begin{aligned}
\dot{u} & =   0, \\  
\dot{y} & =   0, \\
\dot{v} & =  -cv + u^2 w,\\
\dot{w} & =  y + vw - cw, \\
\end{aligned}
\end{equation}
while the nonlinear slow system becomes the so-called {\em reduced problem}
\begin{equation}\label{eq:nonlin-reduced}
\begin{aligned}
u' & =  v, \\
y' & =  - w(1-w), \\
0 & = -cv + u^2w, \\ 
0 & =  y + vw -cw. 
\end{aligned}
\end{equation}
Now we choose appropriate solutions to \cref{eq:nonlin-layer,eq:nonlin-reduced}, and glue them together at their end-states of the dependant variables, producing weak travelling wave solutions to \cref{eq:ardpde} for $\ve =0$. In \cite{Harley_2014}, the authors then exploit GSPT to show that these solutions perturb appropriately in the full nonlinear ODEs given in \cref{eq:twode}. 
\subsection{The layer problem} 
\begin{figure}[ht]
	\centering
	\includegraphics[height=2.5in]{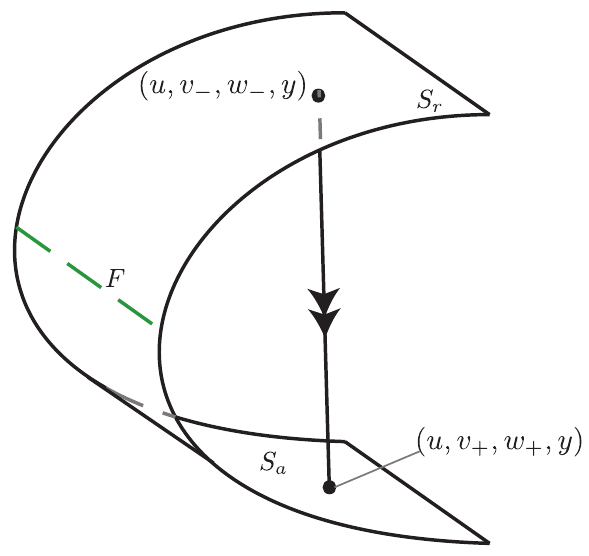}
	\captionsetup{width=.95\linewidth}
	\caption{A schematic of the critical manifold $S$.  The fold curve $F$ is represented by the dashed line (green online).  The upper part of the surface is the repelling side of the manifold $S_r$ and the lower part the attracting side of the manifold $S_a$.  The flow of the layer problem is along {fast fibres}, an example of which is shown.  Fast fibres connect a point on $S_r$ (labelled $(u,v_-,w_-,y)$), to a point of $S_a$ (labelled $(u,v_+,w_+,y)$). Along these fast fibres $u$ and $y$ are constant.  From the layer dynamics, it follows that the direction of the flow can only be from $S_r$ to $S_a$.}
	\label{fig:fast-fibres}\vspace{-0.5cm}
\end{figure}

Steady states of the layer problem given in \cref{eq:nonlin-layer} define a critical manifold $S$, represented as a graph over $(u,w)$,
\begin{equation} \label{eq:critman} 
S = \left\{ (u,v,w,y) \bigg{|}  v = \frac{u^2 w}{c}, y=-\frac{u^2w^2}{c}+cw  \right\},
\end{equation}
and we will henceforth consider the existence problem in a single coordinate chart by projecting onto $(u,w)$ space. The most important property of the critical manifold $S$ is that it is {\em folded}. We cite the following lemma from \cite{Harley_2014} without proof: 
\begin{lemma}[\cite{Harley_2014}, Lem 2.2] The critical manifold $S$ of the layer problem is folded around the curve 
	$$ F(u,w) := 2u^2w - c^2 = 0 $$ 
	in the $(u,w)$ plane with one attracting side $S_a$ and one repelling side $S_r$. 
\end{lemma}
We refer to the curve $F(u,w) =0$ as the {\em fold curve} or the {\em wall of singularities}. The terminology follows from the behaviour of the reduced problem (see below). The so-called {\em fast fibres} of the layer problem connect points on $S$ with constant $u$ and $y$. Due to the stability of $S$, the direction of the flow along these fast fibres is from the repelling side $S_r$ to the attracting side $S_a$ (see \cref{fig:fast-fibres}).

\subsection{The reduced problem} \Cref{eq:nonlin-reduced} is a differential-algebraic problem. The reduced flow is constrained to the critical manifold $S$, and the reduced vector field is contained in the tangent bundle of $S$. Since $S$ is given as a graph over $(u,w)$ space, we study the reduced flow in the single coordinate chart. In \cite{Harley_2014} it was shown that the reduced problem contains a so-called {\em folded saddle canard point} \cite{Wechselberger_Pettet_10}.

Eliminating $v$ and $y$ from \cref{eq:nonlin-reduced} gives  the reduced vector field on $S$, 
\begin{equation}\label{eq:reducedsing}
\begin{pmatrix} c & 0 \\ -2uw^2/c & c - 2u^2w/c \end{pmatrix} \begin{pmatrix} u \\ w \end{pmatrix}' = \begin{pmatrix} u^2w \\ -w(1-w) \end{pmatrix} .
\end{equation}
The left hand side of \cref{eq:reducedsing} is singular along the fold curve $F(u,w) = 0$, but can be desingularised by multiplying both sides by the co-factor matrix of the matrix on the left in \cref{eq:reducedsing}, and by rescaling the independent variable $z = z(\bar{z})$ such that
$$
\frac{dz}{d\bar{z}} = c^2 - 2u^2w.
$$
This gives the {\em desingularised system} 
\begin{equation}\label{eq:desing}
\begin{split}
\frac{du}{d\bar{z}} & = cu^2w - \frac{2u^4w^2}{c} \\ 
\frac{dw}{d\bar{z}} & = -cw(1-w) + \frac{2u^3w^3}{c}. 
\end{split}
\end{equation}

The equilibrium points of \cref{eq:desing} are $(u_U, w_U) = (0,1)$, $(u_S, w_S) = (u_{\infty},0)$, $u_{\infty} \in \mathbb{R}$ and 
\begin{equation}\label{eq:canard-point}
(u_H,w_H) = \left(\frac{c}{4}\left[c + \sqrt{c^2 + 8}\right], \frac{1}{u_H + 1}\right).
\end{equation}
The first two equilibrium points listed correspond to the background states of \cref{eq:ardpde}, while the last is a product of the desingularisation.  More specifically, the Jacobian at $(u_U, w_U) = (0,1)$ has eigenvalues and eigenvectors 
\[ \lambda_1 = c, \quad \boldsymbol{\psi}_1 = (0,1), \quad \lambda_2 = 0, \quad \boldsymbol{\psi}_2 = (1,0), \]
and is therefore centre-unstable; the Jacobian at $(u_S, w_S) = (u_{\infty},0)$ has eigenvalues and eigenvectors
\[ \lambda_1 = -c, \quad \boldsymbol{\psi}_1 = (-u_{\infty}^2,1), \quad \lambda_2 = 0, \quad \boldsymbol{\psi}_2 = (1,0), \]
and is therefore centre-stable; and finally, the Jacobian at $(u_H,w_H)$ has eigenvalues and eigenvectors
\[ \lambda_{\pm} = \left(\frac{c - \sqrt{c^2 + 8}}{2}\right)^4\left[1 \pm c\sqrt{\left(\frac{4}{c - \sqrt{c^2 + 8}}\right)^4 - 3}\right], \quad \boldsymbol{\psi}^{\pm} = (f^{\pm}(c),-1), \]
with 
$$
f^{\pm}(c) := \frac{c^2(c + \Gamma)^4}{64(c^2 + c \Gamma + 1) \pm 2(c + \Gamma)^2\sqrt{16 + 24c \Gamma - 48c^2 + 6c^3\Gamma - 6c^4}},
$$
where $\Gamma: = \sqrt{c^2 + 8}$,
and is therefore a saddle for all $c>0$.   

To obtain the $(u,w)$-phase portrait in terms of the variable $z$, we observe that $\displaystyle \diff{z}{\bar{z}} > 0$ on $S_a$ (that is, below the fold curve $F$), while $\displaystyle \diff{z}{\bar{z}} < 0$ on $S_r$.  Therefore, the direction of the trajectories in the $(u(z),w(z))$-phase portrait will be in the opposite direction to those in the $(u(\bar{z}),w(\bar{z}))$ phase portrait for trajectories on $S_r$, but in the same direction for trajectories on $S_a$. This does not affect the stability or type of the fixed points $(u_U, w_U)$ and $(u_S,w_S)$ as they are on $S_a$.  However, $(u_H,w_H)$ is not a fixed point of \cref{eq:reducedsing}.  Rather, as the direction of the trajectories on $S_r$ are reversed, the saddle equilibrium of \cref{eq:desing} becomes a {folded saddle} canard point of \cref{eq:reducedsing} \cite{Wechselberger_Pettet_10}.  In particular, on $S_r$ the stable (unstable) eigenvector of the saddle equilibrium of \cref{eq:desing} becomes the unstable (stable) eigenvector of the folded saddle canard point.  This allows two trajectories to pass through $(u_H,w_H)$: one from $S_a$ to $S_r$ and one from $S_r$ to $S_a$. The former is the so-called \emph{canard} solution and the latter the \emph{faux canard} solution \cite{Wechselberger_Pettet_10}.  

\begin{figure}[ht]
	\centering
	\includegraphics[height=7cm]{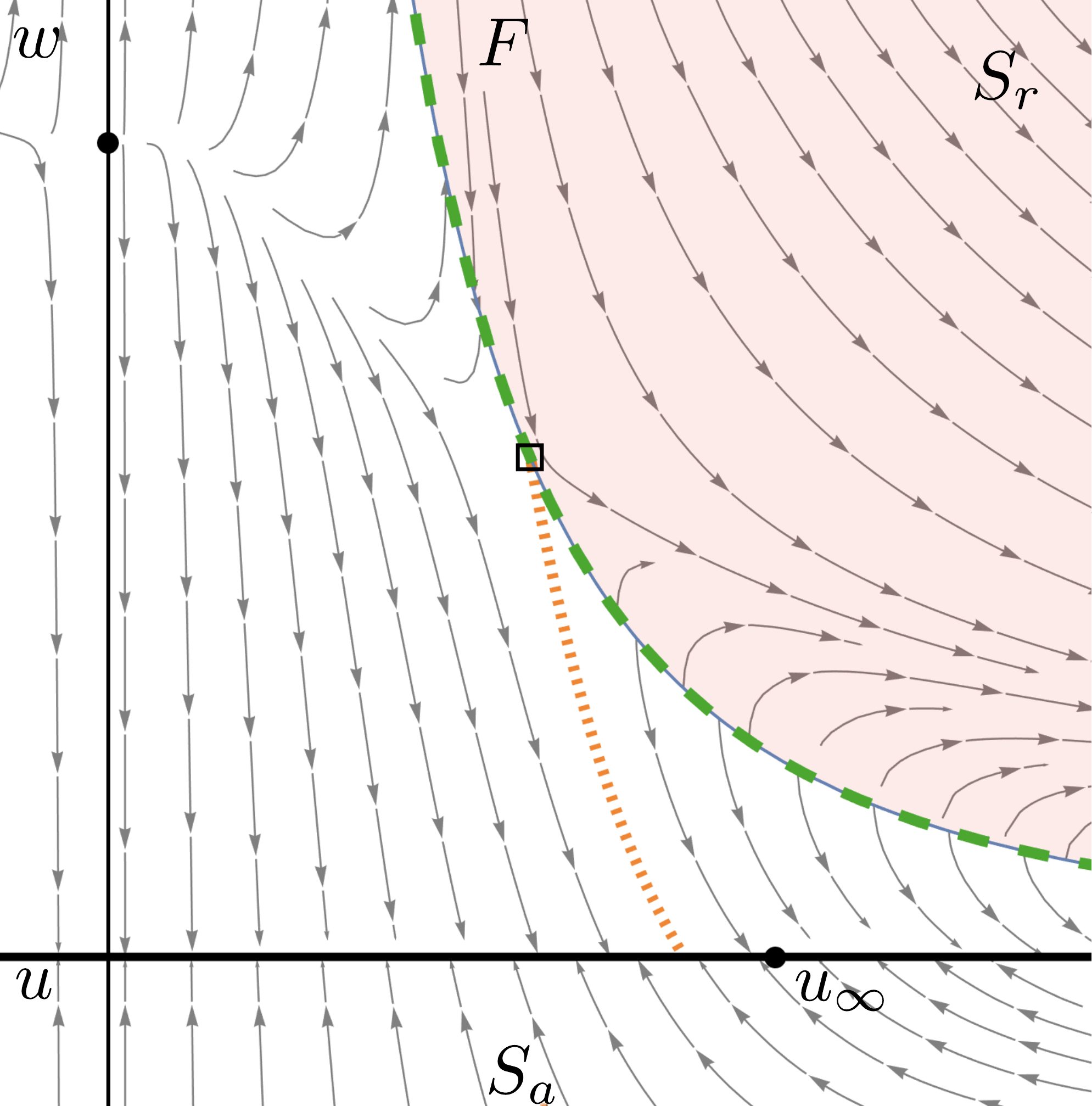}
	\captionsetup{width=.95\linewidth}
	\caption{The $(u,w)$-phase portrait parameterised by the variable $z$. The fold curve (dashed, green online) is labelled $F$ and the folded saddle canard point is the open black square on it. The two solid black circles are the background states $(0,1)$ and $(u_\infty,0)$, which are fixed points of both \cref{eq:reducedsing,eq:desing}. Travelling wave solutions are connections from unstable steady state $(0,1)$ to any of the family of stable steady states $(u_{\infty},0)$ along the $u$-axis. The region below $F$, labelled $S_a$ corresponds to the attracting side of the critical manifold $S$, and above $F$, (red online), corresponds to the repelling side $S_r$.  The dotted line connecting the canard point (orange online) to the line of steady states is a separatrix (faux canard). Thus, existence of a heteroclinic connection (travelling wave) from the left steady state to the point marked $u_\infty$ is only possible if the trajectory passes through the canard point and then travels along the repelling side of the critical manifold before travelling back down to the attracting sheet via a fast fibre. This results in a shock fronted travelling wave.}
	\label{fig:phase-planes}\vspace{-0.75cm}
\end{figure}

The $(u,w)$-phase portrait parameterised by $z$ is shown in \cref{fig:phase-planes}.
\subsection{Travelling wave solutions} 
\begin{figure}
	\makebox[\linewidth][c]{%
			\includegraphics[width=0.8\textwidth]{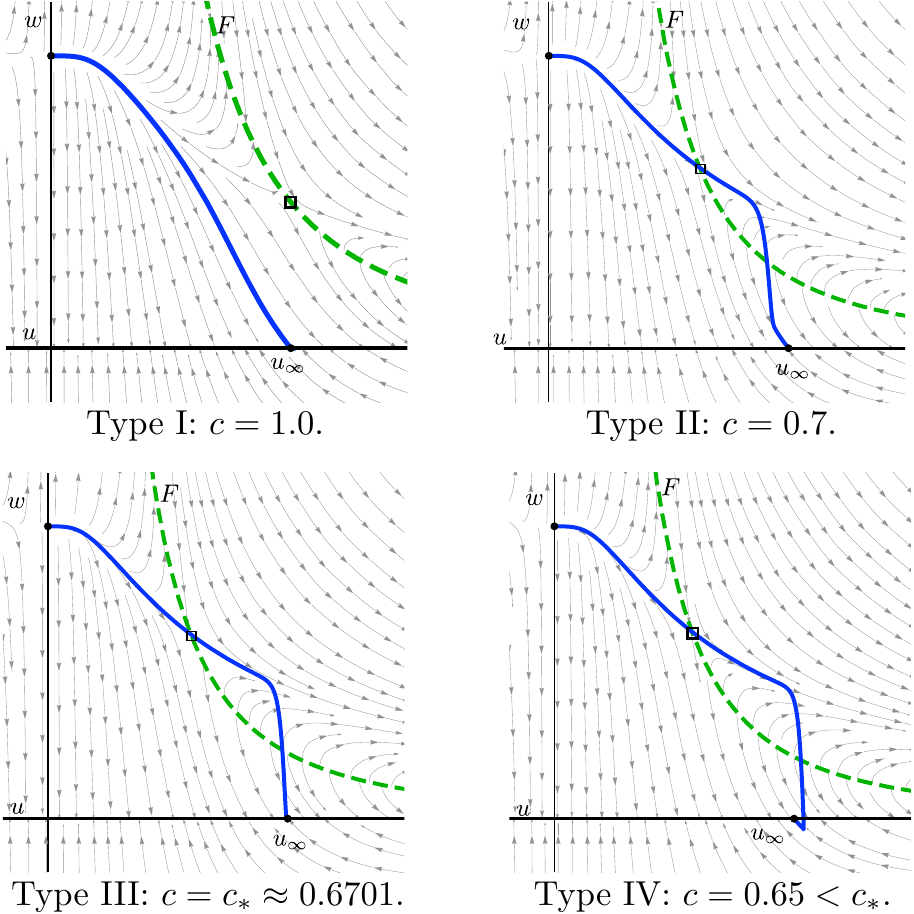}}
				\captionsetup{width=.95\linewidth}
	\caption{An illustration of the four different types of waves found in \cite{Harley_2014} in the phase portrait of the critical manifold $S$ as $c$ is varied, for fixed $u_{\infty} = 1$. The fold lines are indicated by the green dashed lines labelled $F$. As in \cref{fig:phase-planes} the attracting sheet of the critical manifold is to the left of the fold, while the repelling sheet is to the right. Type I waves are smooth and do not cross to the repelling side of $S$. Type II waves are sharp fronted, owing to passing through the canard point on the fold of the critical manifold to the repelling sheet, type IV waves are also sharp-fronted travelling solutions, but are non-monotone. Type III waves, which exist for a unique wave speed $c=c_*$, are the transition between type II and type IV waves where the waves jump through the fast system directly to the line of fixed points on the critical manifold.}\vspace{-.75cm}
	\label{fig:waves}
\end{figure}

As alluded to in the introduction, four distinct types of travelling wave solutions to \cref{eq:ardpde} were identified in \cite{Harley_2014}, denoted types I, II, III, and IV (see \cref{fig:waves223}). The solutions were found as solutions to the desingularised system of the reduced problem and were glued together with (appropriate) fast fibres of the layer problem to produce (weak) traveling wave solutions to the full nonlinear travelling wave PDE given in \cref{eq:twpde} (with $\ve =0$). These solutions were then shown to persist for small enough values of the diffusion parameter $\ve$ via standard approaches in GSPT. \cref{fig:waves} provides an example of the four types of waves found in the phase portrait of their desingularised reduced systems. Type I waves are smooth positive waves lying entirely in the attracting sheet of the critical manifold. Type II waves exhibit a shock in $w$ (in the singular limit). They pass through the folded saddle canard point in the reduced problem, and then travel along a fast fibre of the layer problem, landing on the attracting branch of the critical manifold, from which they continue on to the steady state $u_\infty$. The length of the jump is determined by the wave speed $c$ (or by $u_\infty$) and the symmetry of $S$. In particular the jump in $w$ is symmetric around the fold curve $F$ with $u$ fixed \cite{Harley_2014}. Type III waves are those that jump directly from the repelling sheet of the critical manifold $S$ to the line of steady states of the reduced problem. Type IV waves are those for which $w$ exhibits negative values after the jump.  

\section{The spectral problem, essential and absolute spectrum}\label{sec:spec}
In this section, and what follows, we assume that a travelling wave solution to \cref{eq:ardpde} of type I-IV is given, denoted by $\bu := (u,w)^\top$. We view the travelling wave $\bu$ as a  steady state to \cref{eq:twpde}, and motivated by dynamical systems theory, we want to examine a linear spectral problem associated with \cref{eq:twpde} at $\bu$. The linearisation of \cref{eq:twpde} at $\bu$ is formally given by:
\begin{equation}\label{eq:lin-pde}
\begin{pmatrix} 
p \\ r 
\end{pmatrix}_t 
= \ve \begin{pmatrix} p \\ r \end{pmatrix} '' + c \begin{pmatrix} p \\ r \end{pmatrix}' - \begin{pmatrix} 0 \\  w p'+ u'r  \end{pmatrix}'  + 
\begin{pmatrix} -2uw p -u^2 r \\ (1-2w)r  \end{pmatrix}.
\end{equation}
We denote the linear operator $L(\bu)$ as the right hand side of \cref{eq:lin-pde} acting on the perturbations $p$ and $r$. That is:
$$
L(\bu) := \ve \partial_{zz} + c \partial_z - \begin{pmatrix} 0 & 0 \\ w \partial_{zz} + w' \partial_z & u' \partial _z + u'' \end{pmatrix} + 
\begin{pmatrix}-2uw & -u^2 \\ 0 & (1-2w)  \end{pmatrix}.
$$
We define {\em the spectrum of $L(\bu)$}, denoted $\sigma(L(\bu))$ as those $\lambda \in \C$ such that $L(\bu) - \lambda \I$ is not invertible on the space $\cX := \cH^1(\R) \times \cH^1(\R)$ (that is we require both $p$ and $r$ and their derivatives to be square integrable functions from $\R \to \C$). To find such values of $\lambda$ we study the system of non-autonomous ODEs
\begin{equation}\label{eq:spectralode}
\ve \begin{pmatrix} p \\ r \end{pmatrix} '' + c \begin{pmatrix} p \\ r \end{pmatrix}' - \begin{pmatrix} 0 \\ u'r + w p' \end{pmatrix}'  + 
\begin{pmatrix} (-2uw - \lambda )p -u^2 r \\ (1-2w- \lambda)r  \end{pmatrix} = \begin{pmatrix} 0 \\ 0 \end{pmatrix} 
\end{equation}
The idea now is to use a linearisation of the Li\'enard coordinates introduced in \cref{eq:lienardnonlin} to derive a linear system with the same slow-fast structure as the original travelling waves $\bu$. We introduce the new linearised, Li\'enard variables 
\begin{equation}\label{eq:lienardlin}
q:= p' \textrm{ and }s := \ve r' + cr - u'r - w q,
\end{equation} and we rewrite $\left(L(\bu) - \lambda \I \right) \begin{pmatrix} p \\ r \end{pmatrix} = 0$ as a slow-fast, linear, non-autonomous system with two fast ($q$ and $r$) and two slow ($p$ and $s$) variables
\begin{equation}\label{eq:lin-slow}
\begin{pmatrix} p \\ s \\ \ve q \\ \ve r  \end{pmatrix}' = 
\begin{pmatrix}0 & 0 & 1 & 0 \\ 0 & 0 & 0 & \lambda - 1 + 2w \\ \lambda + 2uw & 0 & -c & u^2 \\ 0 & 1 & w & u'-c   \end{pmatrix}
\begin{pmatrix} p \\ s \\ q \\ r \end{pmatrix}.
\end{equation}
We refer to \cref{eq:lin-slow} as the {\em (linear) slow system}, again with the slow variable $z$.  For notational convenience, we will denote the vector $(p,s,q,r)$ as $\bp$ and note that we can write \cref{eq:lin-slow} as $\bp' = A(z;\lambda,\ve) \bp$ where $A(z;\lambda,\ve)$ is the matrix given by 
\begin{equation}\label{eq:matslow}
A(z;\lambda,\ve) := \begin{pmatrix}0 & 0 & 1 & 0 \\ 0 & 0 & 0 & \lambda - 1 + 2w \\ (\lambda + 2uw)/\ve & 0 & -c/\ve & u^2/\ve \\ 0 & 1/\ve & w/\ve & (u'-c)/\ve   \end{pmatrix}.
\end{equation}
We can make the same change of independent variable as before, $\zeta = z/\ve$, to derive the {\em (linear) fast system}
\begin{equation}\label{eq:lin-fast}
\begin{pmatrix} \dot{p} \\ \dot{ s} \\ \dot{q} \\ \dot{r} 
\end{pmatrix} = \begin{pmatrix} 0 & 0 & \ve & 0 \\ 0 & 0 & 0 & \ve(\lambda - 1 + 2w) \\ \lambda + 2uw & 0 & -c & u^2 \\ 0 & 1 & w & u'-c   \end{pmatrix} \begin{pmatrix} p \\ s \\ q \\ r \end{pmatrix} =: B(\zeta;\lambda,\ve) \bp. 
\end{equation}
We next recall that our travelling waves in both the slow and the fast variables are asymptotically constant - they either satisfy the boundary conditions given in \cref{eq:bcs} or the jump conditions. The jump conditions in this framework are determined by the symmetry of $S$ about the fold curve and are given as 
	\begin{align*}
		v_+ - v_- &= \frac{u^2}{c}(w_+ - w_-), \\
		w_+ + w_-& = \frac{c^2}{u^2} 
	\end{align*}
	where the $\pm$ subscript denotes the value of the given variable at the beginning or end state of the shock respectively and we recall that $u$ is constant during the shock \cite{Harley_2014}.  
	As $z$ or $\zeta \to \pm \infty$ the matrices $A(z;\lambda,\ve)$, and $B(\zeta;\lambda,\ve)$ will tend towards the constant matrices $A_\pm(\lambda,\ve)$ and $B_\pm(\lambda,\ve)$ respectively. The matrices $A_\pm$ are given by: 
	\begin{equation*}
		A_-(\lambda,\ve) :=  \begin{pmatrix}0 & 0 & 1 & 0 \\ 0 & 0 & 0 & \lambda +1 \\ \lambda/\ve & 0 & -c/\ve & 0 \\ 0 & 1/\ve & 1/\ve & -c/\ve   \end{pmatrix}, \, 
		A_+(\lambda,\ve) :=  \begin{pmatrix}0 & 0 & 1 & 0 \\ 0 & 0 & 0 & \lambda - 1 \\ \lambda/\ve & 0 & -c/\ve & u_\infty^2/\ve \\ 0 & 1/\ve & 0 & -c/\ve   \end{pmatrix}.
	\end{equation*}
	The matrices $B_\pm(\lambda,\ve)$ are given by 
	\begin{equation*}
		\begin{array}{ccc}
			B_\pm(\lambda,\ve) & :=  & \begin{pmatrix} 0 & 0 & \ve & 0 \\ 0 & 0 & 0 & \ve(\lambda - 1 + 2w_\pm) \\ 
			\lambda + 2u w_\pm & 0 & -c & u ^2 \\ 
			0 & 1 & w_\pm & v_\pm -c \end{pmatrix}.
		\end{array}
	\end{equation*}
	Where $u$ is a constant in the fast (nonlinear) system, and $v_\pm$ and $w_\pm$ are the jump conditions that must be satisfied along the fast fibres.
	\subsection{Definition of the essential and point spectrum} \mbox{}
	In this section, we follow \cite{Kap_Prom_2013,Sandstede_02}. The spectrum $\sigma(L(\bu))$ splits up 
	into two parts, the {\em point spectrum}, denoted $\sigma_{\textrm{pt}}(L(\bu))$ and the {\em essential spectrum} denoted $\sigma_{\textrm{c}} (L(\bu))$. We define the point spectrum as the values of $\lambda \in \sigma(L(\bu)) $ where $L(\bu) - \lambda$ has a finite dimensional kernel and cokernel, and the {\em index} of $L(\bu)-\lambda$ := dim(kernel) -- dim(cokernel) is zero. We define the essential spectrum as the complement $\sigma_{\textrm{c}}(L(\bu)) := \sigma(L(\bu)) \setminus \sigma_{\textrm{pt}}(L(\bu))$ of the point spectrum.

	The operator $\frac{d}{dz} - A(z;\lambda,\ve)$ is a relatively compact perturbation of the piecewise operator $\frac{d}{dz} - A_\pm(\lambda,\ve)$ for $z \lessgtr 0$ in $\cH^1(\R^\pm)$, (and likewise for the appropriate $B$ matrices). Thus, the essential 
	spectrum is where the Morse indices (dimension of the unstable spatial eigenspace) of the end states are different \cite{Kap_Prom_2013,Sandstede_02}. 
	
\begin{figure}
	\centering
		\includegraphics[scale=0.25]{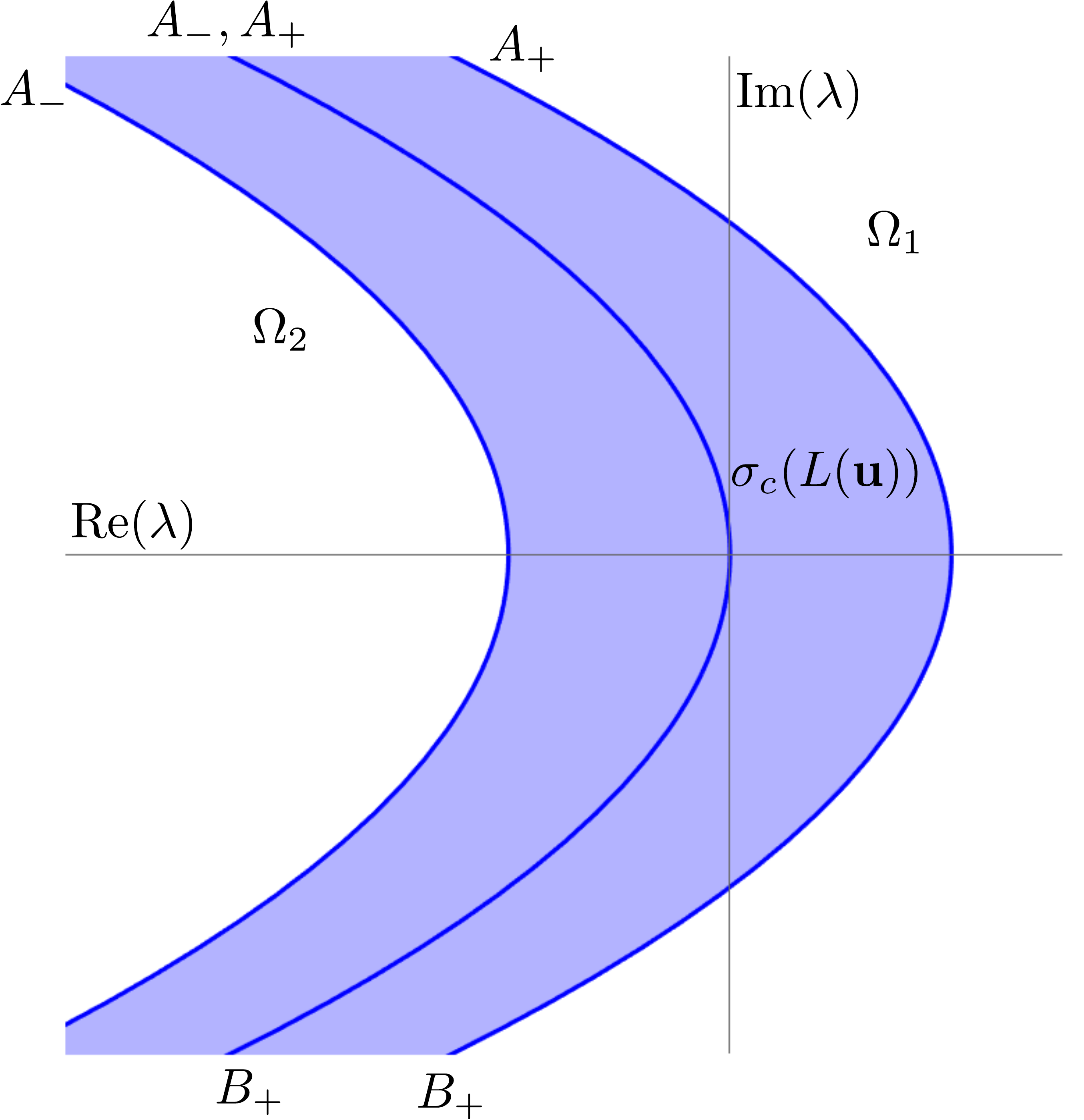}
		\captionsetup{width=.95\linewidth}
		\caption{A plot of the essential spectrum of the operator $L(\bu)$. The dark lines (blue online) bounding the essential spectrum and passing through the origin in the complex plane are the dispersion relations for the matrices $A_\pm$ and $B_+$, labelled accordingly (see \eqref{eq:dispersion}). In all cases qualitatively the essential spectrum is the same. For this figure, the value of $\ve =0.01$ while $c = 1$. The absolute spectrum in this case is the set $(-\infty,1-\frac{c^2}{4\ve}] = (-\infty,-24] \in \R$. In particular it is real, and far to the left (in the region $\Omega_2$ and out of the figure).}\label{fig:essential}
		\vspace{-0.5cm}
	\end{figure}

For waves of type I, II, and IV the end-states of the wave are in the slow system, and so the matrices $A_\pm(\lambda,\ve)$ determine the essential spectrum.  We have that $\lambda \in \sigma_\textrm{c}(L(\bu))$ when $A_+(\lambda,\ve)$ has a different number of unstable spatial eigenvalues from $A_-(\lambda,\ve)$, or either one has a purely imaginary eigenvalue. In all cases, this is a region in the complex plane bounded by the so-called {\em dispersion relations}. These are curves where $A_+(\lambda,\ve)$, $A_-(\lambda,\ve)$ have purely imaginary eigenvalues $ik$ for $k \in \R$, and are the following four curves (two lie on top of each other):  
	\begin{equation}\label{eq:dispersion}
	\begin{array}{lr}
	\lambda  = -\ve k^2 - 1  + i c k , & (A_-(\lambda,\ve) \textrm{ has eigenvalue $ik$}) \\ 
	\lambda  = - \ve k^2 + i c k, &  (A_\pm(\lambda,\ve) \textrm{ has eigenvalue $ik$})\\
	\lambda  = 1-\ve k^2 + i ck & (A_+(\lambda,\ve) \textrm{ has eigenvalue $ik$}) \\
	\end{array}
	\end{equation}
For waves of type III, the end-state of the wave is in the slow system as $z \to -\infty$ but in the fast system as $\zeta \to + \infty$, and now the essential spectrum is the $\lambda \in \C$ when $A_-(\lambda, \ve)$ has a different number of unstable eigenvalues from $B_+(\lambda,\ve)$. 
We note that it is not strictly necessary to use $B_+$ in order to apply Weyl's theorem to compute the essential spectrum of the type III waves, as long as $\ve >0$, due to the equivalence of the fast and slow systems. Indeed, it turns out that the dispersion relations from the matrix $B_+$ for $\ve >0$ define the same set of curves in the spectral parameter as those from $A_+$. This is reflected in the specific values that the jump conditions take for the type III waves ($v_+ = w_+ = 0$).  The dispersion relations for the type III waves are 
\begin{equation}\label{eq:dispersionIII}
	\begin{array}{lr}
	\lambda  = -\ve k^2 - 1  + i c k , & (A_-(\lambda,\ve) \textrm{ has eigenvalue $ik$}) \\ 
	\lambda  = - \ve k^2 + i c k, &  (A_-(\lambda,\ve) \textrm{ has eigenvalue $ik$})\\
	\ve \lambda  = - k^2 + i c k & (B_+(\lambda,\ve) \textrm{ has eigenvalue $ik$}) \\ 
	\ve \lambda  = \ve- k^2 + i ck &  (B_+(\lambda,\ve) \textrm{ has eigenvalue $ik$}). \\
	\end{array}
	\end{equation}
The second and third curves lie on top of each other, even though their expressions are different. The essential spectrum for a type III waves is thus the same as that of types I, II and IV (see \cref{fig:essential}). 
	
We also remark that the dispersion relations divide the complex plane into three disjoint regions. The first we denote by $\Omega_1$. In the type I, II or IV case, this is the region where if $\imag{\lambda} = c k $ for some $k \in \R$, then $\real{\lambda} > 1-\ve k^2,$ i.e. to the right of the essential spectrum.  $\Omega_1$ is also to the right of the essential spectrum in the type III case, though here if $\imag{\lambda} = \frac{ck}{\ve}$, then we require $\real{\lambda}>1- \frac{k^2}{\ve}$. The next region is $\sige{L(\bu)}$ where $L(\bu)- \lambda$ does not have Fredholm index 0. The third remaining region of the complex plane, to the left of $\sige{L(\bu)}$, we denote $\Omega_2$ (see \cref{fig:essential}).
	
	Since we are concerned with stability of the travelling waves found in \cite{Harley_2014}, it is worth mentioning that for all types of travelling waves identified, the intersection of the 
	essential spectrum with the right half plane is nonempty. However, by considering appropriate {\em weights} and {\em weighted spaces} we can move the spectrum of the linearised 
	operator into the left half plane for all four types of travelling waves. 
	For a given weight function, $\tilde{\alpha}(x)$, we define 
	$$
	|| f || _{\cH_{\tilde{\alpha}}^k} := || \tilde{\alpha} f ||_{\cH^k}.
	$$
	For the travelling waves at hand, the essential spectrum due to $A_-(\lambda, \ve)$ is contained in the left half plane, while for the essential spectrum coming from $A_+(\lambda, \ve)$ or $B_+(\lambda,\ve)$, determination of the appropriate weighted space is identical to determining the appropriately weighted space for travelling waves in Fisher's equation. Consequently the appropriate weighted space for travelling waves of all types is given by a so-called two-sided weight
	$$
	\alpha(x) := \left\{ \begin{split} 1 & \quad \textrm{if} \quad z \leq 0 \\ e^{\nu x} & \quad \textrm{if} \quad z >0 \end{split} \right. 
	$$
	with 
	$$
	\nu \in \left( \frac{c - \sqrt{c^2-4\ve}}{2 \ve}, \frac{c + \sqrt{c^2-4 \ve}}{2 \ve} \right).
	$$
	Thus if $\bp \in \cH^1_\alpha$, we have that the essential spectrum of $\frac{d}{dz} - A(z;\lambda,\ve)$ will be contained in the left half plane.

	This implies the presence of a so-called {\em transient,} or {\em convective instability,} \cite{Sandstede_02,sandscheel00}  where small perturbations 
	either outrun the travelling wave, or die back into the wave, resulting in temporal evolution to a translate (perhaps with a slightly modified wave speed) of the original wave. As the perturbation outruns the wave, it can (and generically will) affect the asymptotic decay ratef which, because this equation shares dynamical qualitative (and quantitative) features with Fisher`s equation, will affect the asymptotic wave speed and the position of the centre of the wave, see also \cite{Harley_2014}. The effect is that
	small perturbations of the original travelling wave evolve into waves that are similar in appearance and behaviour to the original wave (even if the difference in an $\cH^1$ norm grows in time), and so we do not really consider these to be instabilities. 	
	What does pose a problem for (spectral) stability is the so-called {\em{absolute spectrum}}. The absolute spectrum is not spectrum per se, but rather is defined as the values of the spectral parameter $\lambda$ where 
	a pair of eigenvalues of the limiting matrices, (i.e. $A_\pm(\lambda,\ve)$ in the type I, II and IV cases and $A_-(\lambda,\ve)$ and $B_+(\lambda,\ve)$ in the type III case) have equal real parts. The absolute spectrum provides a bound for how far the essential spectrum can be moved by considering perturbations with different weights. In particular if the absolute spectrum is in the right half of the complex plane, there is no choice of a weight that can move the essential spectrum into the left half plane. 
	
	The eigenvalues of $A_-(\lambda,\ve)$ for all types of waves are found to be the following,
	\begin{equation}\label{eq:defmu}
	\mu_0^{\pm} := \frac{-c \pm \sqrt{c^2 + 4 \ve \lambda}}{2 \ve} \quad \mu_{-1}^{\pm} :=  \frac{-c \pm \sqrt{c^2+ 4 \ve (\lambda +1)} }{2 \ve}, 
	\end{equation}
	while the eigenvalues of $A_+(\lambda,\ve)$ (for types I, II and IV only) are
	\begin{equation}\label{eq:defrho}
	\rho_0^{\pm} := \mu_0^{\pm} =  \frac{-c \pm \sqrt{c^2 + 4 \ve \lambda}}{2 \ve} \quad \rho_1^{\pm} :=  \mu_1^{\pm} :=  \frac{-c \pm \sqrt{c^2+ 4 \ve (\lambda -1)} }{2 \ve},
	\end{equation}
	and the eigenvalues of $B_+(\lambda,\ve)$ for a type III wave are
	\begin{equation}\label{eq:defbeta}
	\beta_0^{\pm} := \ve \mu_0^{\pm}=  \frac{-c \pm \sqrt{c^2 + 4 \ve \lambda}}{2} \quad \beta_1^{\pm} := \ve \mu_1^{\pm} =  \frac{-c \pm \sqrt{c^2+ 4 \ve (\lambda -1)} }{2}.
	\end{equation}
	The naming conventions are as follows: $\mu$ for $A$ at minus infinity, $\rho$ for $A$ at plus infinity, and $\beta$ for $B$ at plus infinity. The $\pm$ refers to the choice of the square root in the eigenvalue calculation, and the subscript $-1,1,0$ refers to the value of $\lambda$ which makes the eigenvalue with the positive square root $=0$. 
	
	The absolute spectrum is real for all waves and consists of the half line
	\begin{equation}
	\label{eq:absspecdef}
	\sigma_{\textrm{abs}} := \left(-\infty, 1 - \dfrac{c^2}{4\ve} \right],
	\end{equation}
	and hence will be in the left half of the complex plane provided that $c^2>4\ve$. This is identical to the case of the travelling waves found in the Fisher-KPP waves (where $\ve$ is the diffusion parameter/coefficient). However, unlike in the Fisher-KPP case where the diffusion coefficient is often taken to be on the same order as the wave speed, here we have that $0<\ve\ll1$ and so for the parameter regime considered in this manuscript we do not expect the absolute spectrum to destabilise the travelling waves of interest. In the travelling waves of type I-IV studied here, as we shall see, there is another destabilising factor due to an element of the point spectrum entering into the right half plane.

	\section{Point spectrum and the Riccati-Evans function}\label{sec:point}
	
	We next compute the point spectrum, or lack thereof, in the right half complex plane of the linearised operator associated with the travelling waves of types I-IV found in \cref{sec:spec}. To do this, we use a modified version of the so-called {\em Evans function} \cite{Kap_Prom_2013}. In order to verify the lack of point spectrum of travelling waves of type I-III in the right half plane, and to show the existence of an eigenvalue in the case of a type IV wave, we want to exploit the geometry of the system in order to more efficiently make the computations. This results in relating the Evans function to the so-called Riccati equation on the Grassmannian of two planes in $\C^4$. We produce an Evans function of sorts in that it is an eigenvalue detector, though it does not have all the nice properties of the classical Evans function. In particular it is meromorphic rather than analytic, and it does not appear to be independent of the value of $z$ at which it is evaluated. However we show that the zeros of this function are indeed independent of the point of evaluation and provided certain conditions are met, coincide with the multiplicity of the zeros of the Evans function. 
		
	We recall some familiar results arising in the definition of the Evans function that will be useful for our purposes later. For a detailed discussion and proofs, see \cite{Kap_Prom_2013}. We begin with point spectrum that is away from the essential spectrum. We say that $\lambda \not \in \sigma_c(L(\bu))$ is an {\em eigenvalue of the wave $\bu$ (or of $L(\bu)$) } if we can find functions $\begin{pmatrix} \phi_1 \\ \phi_2 \end{pmatrix} \in \cX$ such that $L(\bu)\begin{pmatrix} \phi_1 \\ \phi_2 \end{pmatrix} = \lambda \begin{pmatrix} \phi_1 \\ \phi_2 \end{pmatrix}$.  
	For $\ve \neq 0$, this is equivalent to finding a $\lambda$ for which there is a solution to the linearised slow problem (i.e. a solution to \cref{eq:lin-slow} in the case of a type I wave), or slow--fast--slow problem (a solution to \cref{eq:lin-slow}, then \cref{eq:lin-fast} and then \cref{eq:lin-slow} in the type II and IV case) or slow--fast problem (a solution to \cref{eq:lin-slow}, then \cref{eq:lin-fast} in the type III case) decaying to zero as $z\to \pm \infty$. Exponential dichotomy for $\lambda \in \Omega_1$ means that there is only one way to do this. Let $\Xi^u$ denote the unstable subspace of $A_-(\lambda,\ve)$ and $\Xi^s$ denote the stable subspace of $A_+(\lambda,\ve)$ in the case that $\bu$ is a type I, II, or IV wave, or the stable subspace of $B_+(\lambda,\ve)$ in the case of a type III wave.

	We note that $\Xi^{u,s}$ are each two-dimensional for $\lambda \in \Omega_1$ (to the right of the essential spectrum) while for $\lambda \in \Omega_2$ (to the left fo the essential spectrum) $\Xi^u$ is zero. We thus (initially) restrict our search for eigenvalues to those $\lambda \in \Omega_1$ which are to the right of the essential spectrum. That is, for a $\lambda \in \Omega_1$, we let $W^{u,s}(z)$ be the (two dimensional) span of solutions to the linearised system along a travelling wave decaying to $\Xi^{u,s}$ respectively (the span of the Jost solutions as in \cite{Kap_Prom_2013}).  We have the following: 
		
	\begin{lemma}[\cite{Kap_Prom_2013}] \label{lem:eigen} Let $\lambda \in \Omega_1$, then $W^u(z_0) \cap W^s(z_0) \neq \left\{ 0 \right\} $ for all $z_0 \in \R$ if and only if $\lambda$ is an eigenvalue. 
		\end{lemma}
		
	Now suppose we pick a pair of linearly independent solutions in each of $W^u$ and $W^s$ respectively, then the above lemma says that if we evaluate them at a given fixed $z_0$ (say $z_0 = 0$), then $\lambda$ will be an eigenvalue if and only if the four are linearly dependent. Denoting these solutions by  $\bx^u_1(z;\lambda), \bx^u_2(z;\lambda), \bx^s_1(z,\lambda)$ and $\bx^s_2(z;\lambda)$ We define the {\em Evans function} as 
	\begin{equation}\label{eq:defevans}
		D(\lambda) := \det \begin{pmatrix} \bx^u_1(0;\lambda), \bx^u_2(0;\lambda), \bx^s_1(0,\lambda),\bx^s_2(0;\lambda) \end{pmatrix}
	\end{equation} 
		
	We have the following 
	\begin{theorem}[\cite{Kap_Prom_2013}]
		The functions $\bx^{u,s}_{1,2}(z)$ can be chosen so that $D(\lambda)$ is analytic for $\lambda$ away from the essential spectrum. The roots of the Evans function $D(\lambda)$ are independent of the choice of $z_0$ being chosen to be $0$. The Evans function is unique up to multiplication by a nonzero function $g(\lambda)$. For $\lambda$ to the right of the essential spectrum, the Evans function is zero if and only if $\lambda$ is an eigenvalue of $\bu$. 
	\end{theorem}
		
	We remark that the additional exponential factor present in many Evans function computations \cite{Kap_Prom_2013} is dropped, as in \cite{ledoux2010grassmannian} as the evolution on the Grassmannian will make it redundant. 
		
	\subsection{The Riccati equation and the Grassmannian}In this section, for the description of the Riccati flow on the Grassmanian, we mostly follow, \cite{hermann1977applications,ledoux2010grassmannian,Shayman86} with some small adaptations to make things more clear for our purposes. We want to exploit some of the geometry behind linear ODEs \cref{eq:lin-slow,eq:lin-fast}. The first observation is that because our ODE is linear, the solution operator maps subspaces to subspaces. This means that for $\lambda$ to the right of the essential spectrum, both $W^u(z)$ and $W^{s}(z)$ will each be two dimensional subspaces of $\C^4$ for all $z \in \R$. Since we are interested in tracking the evolution of the entire subspace, we can consider the (nonlinear) ODE on the space of complex two dimensional subspaces of $\C^4$, the Grassmannian of two planes in four space,\cite{hermann1977applications} which we denote $\mathds{G}\textrm{r}(2,4)$. In this manuscript, since we are primarily only considering the Grassmannian of two planes in four space we  drop the numbers and refer to it just as $\G$. Before we describe the associated Riccati equation on $\G$, we pause for a moment to recall some facts about $\G$ and its coordinatisation. These facts (or equivalent generalisations) can be found in most introductory texts on algebraic geometry, see for example \cite{Harris_92,Shaf_94}.
		
	The manifold $\G$ is a smooth, compact, complex manifold, of complex dimension $4$. It is a homogeneous space, $\G \approx \U(4) / (\U(2) \times \U(2))$, where $\U(n)$ is the unitary group - the real Lie group of real dimension $n^2$ of complex matrices $U$ such that $\bar{U}^TU= \I$. We construct charts on the Grassmannian in the usual way, via the Pl\"ucker coordinates. For a pair of vectors $\bv = (v_1, v_2,v_3,v_4)^\top$ and $\bw = (w_1,w_2,w_3,w_4)^\top$, in $\C^4$ we observe that $\bv$ and $\bw$ are linearly independent (i.e. the plane $P_{\bv,\bw}$ spanned by $\bv$ and $\bw$ is an element of $\G$), if and only if the values of $\bK_{ij} := v_iw_j-v_jw_i$ are not all zero for all $i \neq j$. That is the vector $\left( \bK_{12},\bK_{13},\bK_{14},\bK_{23},\bK_{24},\bK_{34} \right) \neq 0$. This naturally embeds $\G$ into $\PP^5$, the complex projective space (this is called the {\em Pl\"ucker embedding}). We will use the usual designation of coordinates in projective space, $\left[ \bK_{12}:\bK_{13}:\bK_{14}:\bK_{23}:\bK_{24}:\bK_{34} \right] $
	to signify that they are not all zero. It can be checked that if $P_{\bv,\bw}$ represents a complex two plane in four space, then the following {\em Pl\"ucker relation} must hold in the Pl\"ucker coordinates: $\bK_{12}\bK_{34} - \bK_{13}\bK_{24} + \bK_{14}\bK_{23} = 0$. In this way, $\G$ is seen to be a smooth (because it is a homogeneous space) variety in $\PP^5$ of complex projective space. This also gives it the structure of a complex manifold. In a given chart, we can view $\G$ as a graph over the remaining variables. For example, suppose that $\bK_{12} \neq 0$,  then in the Pl\"ucker coordinates we have, by dividing through by $\bK_{12}$ , that our plane is represented by the sextuplet $[1:\bK_{13}:\bK_{14}:\bK_{23}:\bK_{24}:\bK_{13}\bK_{24}-\bK_{14}\bK_{23}],$ and that this represents the plane spanned by $(1,0,-\bK_{23}, -\bK_{24})^\top$ and $(0,1,\bK_{13},\bK_{14})^\top$, which we will write in so-called {\em frame notation}  \cite{ledoux2010grassmannian,Shayman86}
		$$
		\begin{pmatrix} 1 & 0 \\ 0 & 1 \\ -\bK_{23} & \bK_{13} \\ -\bK_{24} & \bK_{14} \end{pmatrix} = \begin{pmatrix} \I \\ \bK \end{pmatrix}. 
		$$
		
	The $4 \times 2$ matrix written as a pair of $2 \times 2$ matrices is called a {\em frame} for the plane that is the span of its columns. Now we want to see how our linear ODE induces a flow on $\G$. Such a flow will be called {\em the associated Riccati equation}. We describe the general process, and then later consider the linear equation coming from the spectral problem at hand. We begin by considering a $4 \times 4$ linear ODE acting on {\em pairs} of vector spaces, and writing it in the frame notation form that will be useful later \cite{hermann1977applications,ledoux2010grassmannian,Shayman86}: 
	\begin{equation}\label{eq:seventeen}
		\begin{bmatrix} \bX \\ \bY \end{bmatrix}' = \A(z)  \begin{bmatrix} \bX \\ \bY \end{bmatrix}  := \begin{bmatrix} A(z) & B(z)  \\ C(z) & D(z) \end{bmatrix}\begin{bmatrix} \bX \\ \bY \end{bmatrix}
	\end{equation}
	where $\bX, \bY, A, B, C, D$ are all $2 \times 2$ matrices in the independent variable $z$.  
		
	Suppose, for the moment that our evolution takes place where $\bX(z)$ is invertible. We can therefore represent the plane $\begin{bmatrix} \bX \\ \bY \end{bmatrix}$ by the plane $\begin{bmatrix} \Id \\ \bY \bX^{-1} \end{bmatrix}$. Denoting the matrix $\bY \bX^{-1}$ by $\bW$, we have that 
	\begin{equation}
		\begin{split}
		\bW' & = \left( \bY \bX^{-1} \right)' = \bY'\bX^{-1} + \bY(\bX^{-1})' \\ 
		& = \bY'\bX^{-1} - \bY\bX^{-1}\bX'\bX^{-1}  \\ 
		& = \left( C \bX + D \bY\right) \bX^{-1} - \bY\bX^{-1}\left(A\bX + B\bY \right) \bX^{-1} \\ 
		\end{split}
	\end{equation}
	where the second step used the fact that $\bX \bX^{-1} = \I$ and the third used \cref{eq:seventeen}. 
	Substituting back in gives 
	\begin{equation} \label{eq:generic}
		\bW'  =  C + D \bW - \bW A  - \bW B\bW.
	\end{equation}
	\Cref{eq:generic} will be called {\em the (associated) Riccati equation} \cite{levin1959matrix,martin1978applications,Shayman86}. It is a higher order analogue of the familiar Riccati equation for second order linear ODEs. This Riccati equation is a nonlinear, non-autonomous ODE of half of the original order. 
	The Riccati equation as written in \cref{eq:generic} governs the flow on a chart of $\G$ equivalent to the original flow prescribed by \cref{eq:seventeen}.  Just as in the more familiar lower order case, solutions to the Riccati equation can become infinite \cite{levin1959matrix}. 
	Geometrically, this means that we are leaving the chart of $\G$ (as $\det(\bX) \to 0$) \cite{ledoux2010grassmannian}. We will return to how to handle this later, but for the moment, we wish to understand how the Evans function defined above fits into the Riccati equation formulation.  
		
	The spans of solutions $W^{u,s}(z)$ decaying to $\Xi^{u,s}$ as $z \to \pm \infty$ are solutions to the Riccati flow on $\G$.  We write them as $ \begin{bmatrix}\bX^u \\ \bY^u \end{bmatrix}$, for the span of $W^u(z)$ and  $ \begin{bmatrix}\bX^s \\ \bY^s \end{bmatrix}$  for the span of $W^s(z)$ where $\bX^{u,s}$ and $\bY^{u,s}$ are each $2 \times 2$ matrices (the pair $\bX^{u,s}$ and $\bY^{u,s}$ are called the Jost matrices in \cite{Kap_Prom_2013}), and again, assuming that we stay in the same chart (i.e $\det (\bX^{u,s}) \neq 0$), we have two solutions to the Riccati flow, $\bW^u(z) := \bY^u(\bX^u)^{-1}$ and $\bW^{s}(z):= \bY^s(\bX^s)^{-1}$. Recall that the eigenvalue problem as we have set it up is to determine whether or not the subspaces $W^{u,s}(z_0)$ intersect nontrivially. So writing the definition of the Evan's function from \cref{eq:defevans} in this new notation, we are interested in zeros of the following function:
	$$D(\lambda) := \det \begin{bmatrix} \bX^u(z_0,\lambda) & \bX^s (z_0,\lambda)\\ \bY^u(z_0,\lambda) & \bY^s(z_0,\lambda) \end{bmatrix}, $$ and we know that the subspaces represented by $\begin{bmatrix} \bX^{u,s}(z_0,\lambda) \\ \bY^{u,s}(z_0,\lambda)\end{bmatrix} $ are the same as those represented by $\begin{bmatrix} \Id \\ \bW^{u,s}(z_0,\lambda) \end{bmatrix}$. The question is how to relate the determinant of $\begin{bmatrix} \Id & \Id \\ \bW^u(z_0,\lambda) & \bW^s(z_0,\lambda)\end{bmatrix}$ to $D(\lambda)$?
		
	It is straightforward to check that for a pair of $2 \times 2$ matrices $A$ and $B,$ the following holds
	\begin{equation}\label{eq:detthm}
		\det \begin{pmatrix} \Id & \Id \\ A & B \end{pmatrix} = \det(B- A). 
	\end{equation}
	That is, the determinant of the $4 \times 4 $ matrix on the left is equal to the determinant of the difference of the matrices $B$ and $A$. 
	This is in fact generically true for $n \times n$ matrices, one just replaces the $2 \times 2$ with the appropriately sized identity matrix. It can also be extended to matrices with a block structure of a more generic type (see \cite{Silvester_2000}), though we will not need the full generic statement here. We thus have: 
	$$ \det \begin{bmatrix} \Id & \Id \\ \bW^u(z_0,\lambda) & \bW^s(z_0,\lambda)\end{bmatrix} = \det(\bW^s(z_0,\lambda) - \bW^u(z_0,\lambda)).$$
	Denote the function
	\begin{equation}\label{eq:defe}
		E(z_0;\lambda) := \det(\bW^s(z_0;\lambda) - \bW^u(z_0;\lambda)).
	\end{equation} 
	Next, we note that  
	$$\begin{bmatrix} \Id & \Id \\ \bW^u(z_0,\lambda) & \bW^s(z_0,\lambda)\end{bmatrix} =  \begin{bmatrix} \bX^u(z_0,\lambda) & \bX^s(z_0,\lambda)\\ \bY^u(z_0,\lambda) & \bY^s(z_0,\lambda) \end{bmatrix} \begin{bmatrix} 
	(\bX^u)^{-1}(z_0,\lambda) & 0 \\  0 & (\bX^s)^{-1}(z_0,\lambda) \end{bmatrix} $$ and taking determinants and using \eqref{eq:detthm} we have that 
	$$ \det(\bX^u(z_0;\lambda)) \det(\bX^s(z_0;\lambda)) E(z_0;\lambda) = D(\lambda).$$
	\begin{definition}
		We call the function $E(z_0,\lambda)$ {\em the Riccati-Evans function.} 
	\end{definition}

		\subsection{Changing charts} In this section, we use the general coordinatisaion of the Grassmannian found in \cite{Shaf_94}. A chart on the Grassmannian is a map $T: \G \to \C^4$. We can think of the charts as parametrised by invertible matrices $\bT \in GL(\C,4)$ in the sense that if we multiply a frame $\begin{pmatrix} \bX \\ \bY \end{pmatrix}$ by a matrix $\bT$ and then compose the result with the Pl\"ucker coordinate map, we get a new coordinate representation for the original plane. For example, suppose we consider the plane spanned by the columns of the frame $\begin{pmatrix} 0 \\ \I  \end{pmatrix}$. This plane is not in the chart where $\bK_{12} \neq 0$ described earlier, rather its coordinates in $\PP^5$ are $[0:0:0:0:0:1]$, so it lies in the chart where $\bK_{34} \neq 0$. However if we multiply the original frame by the matrix $\bT = \begin{pmatrix} 0 & \I \\ \I & 0 \end{pmatrix}$, then in the new coordinate chart associated with $\bT$ we have that the frame is given as $\begin{pmatrix} \I \\ 0 \end{pmatrix}$, and so in this chart, the same plane is represented by $\bK_{12} \neq 0$. This parametrisation has several advantages, namely it allows us to write down a single expression for the evolution of an ODE which changes implicitly depending on the chart (matrix $\bT$) we choose.  
		
		We next write out our matrix Riccati equation in the chart parametrised by $\bT$. This is the evolution equation on $\G$ under the change of variables determined by $\bT$. Suppose that in our original variables 
		\begin{equation}
		\begin{bmatrix} \bX \\ \bY \end{bmatrix}' = \A(z)  \begin{bmatrix} \bX \\ \bY \end{bmatrix}
		\end{equation} 
		Then if $\bT$ is an invertible matrix, so that we have 
		$$ \bT \begin{bmatrix} \bX \\ \bY \end{bmatrix} =: \begin{bmatrix} \bX_\bT \\ \bY_\bT \end{bmatrix} $$ 
		and 
		\begin{equation}
		\begin{bmatrix}  \bX_\bT \\ \bY_\bT \end{bmatrix}' = \bT \A(z) \bT^{-1}  \begin{bmatrix}  \bX_\bT \\ \bY_\bT \end{bmatrix} =: \begin{bmatrix} A_\bT(z) & B_\bT(z) \\ C_\bT(z) & D_\bT(z) \end{bmatrix}\begin{bmatrix} \bX_\bT \\ \bY_\bT \end{bmatrix}.
		\end{equation} 
		Defining $\bW_T = \bY_\bT \bX_\bT^{-1}$, the Riccati equation in this chart is 
		$$ 
		\bW_{\bT}'  =  C_\bT + D_\bT \bW_{\bT} - \bW_{\bT} A_\bT  - \bW_{\bT} B_\bT\bW_{\bT}.
		$$
		We have therefore absorbed the chart implicitly into the computations, in order to have a single set of ODEs to evolve.
		
		Likewise, we can define the Riccati-Evans function on this chart 
		$$
		E_\bT(z_0;\lambda) := \det(\bW_\bT^s(z_0;\lambda) - \bW_\bT^u(z_0;\lambda)),
		$$
		and the relation
		\begin{equation}\label{eq:ric-evans-evans}
		\det(\bT^{-1}) \det(\bX_\bT^u(z_0;\lambda)) \det(\bX_\bT^s(z_0;\lambda)) E_\bT(z_0;\lambda) = D(\lambda)
		\end{equation}
		still holds. The Riccati-Evans function is not independent of the change of coordinates, but we use this to our advantage. We will choose a chart (matrix $\bT$) so that $\det(\bT) =1$ and $\det(\bX_\bT^{u,s}) \neq 0$ in the spectral parameter regime of interest, and produce a function $E_\bT$, the zeros of which coincide with those of $D(\lambda)$. 
		
		We note that in the current notation, the function defined in \cref{eq:defe} is for the chart corresponding to the identity. That is
		$$ 
		E(z_0;\lambda) = E_\I(z_0;\lambda).
		$$

		\subsection{Extension into the essential spectrum}
		Using $\Xi^{u,s}$ defined above as initial conditions, we can then (numerically) compute the Riccati-Evans function on any chart associated with an invertible matrix $\bT$ for any $\lambda \in \Omega_1$. We would like to consider a larger domain of $\lambda \in \C$ however, not just those $\lambda \in \Omega_1$. This is relatively straightforward provided we stay away from values of $\lambda$ in the absolute spectrum, computed above in \cref{eq:absspecdef}.
		
		To extend the Evans function, we track the eigenvectors associated with $\mu_{0,-1}^+$ and $\rho_{0,1}^-$ (see \cref{eq:defmu,eq:defrho}) as we vary $\lambda$. Starting with a $\lambda \in \Omega_1$, we can continue the Evans function (and the Riccati-Evans function) as we vary $\lambda$ through the curves defined by the dispersion relations in \cref{eq:dispersion}. A root of $D(\lambda)$ will no longer be evidence of {\em any} solution which decays at $\pm \infty$ but rather a solution that decays at $\pm \infty$ along the eigenspaces $\Xi^{s,u}$. For example, the eigenvalue associated with the derivative of the type I, II and IV waves found in \cref{sec:back} will not be a root of this extended Riccati-Evans function, as the solution will not decay along the appropriate subspace. So, even though $\lambda = 0$ (and in fact any $\lambda \in \sige{L}$ not on the boundary of $\sige{L}$) will technically be an eigenvalue of $L$, in the sense that there will be a decaying $L^2$ solution to the ODE, it will not be a root of this extended Evans function. In some sense this is preferred as roots of the Evans function found in this manner can not be removed by considering functions in weighted space which moves the essential spectrum into the left half plane, whereas eigenvalues which are removed due to weighting are associated with so-called {transient or convective} instabilities \cite{Kap_Prom_2013,sandscheel00} which are known to affect the temporal dynamics of the wave less strongly or noticeably than eigenvalues which cannot be weighted away. As we shall see, it is these roots of the extended Evans function which are associated with a change in stability of the travelling waves outlined in \cref{sec:back}. 
		
		\subsection{Winding numbers} One typical way that the analyticity of the Evans function $D(\lambda)$ is employed is via the argument principle from complex analysis. This can be stated as follows 
		\begin{theorem}[\cite{carrier2005functions}]
			Suppose $f:\Omega \to \C$ is a complex meromorphic function on a simply connected domain $\Omega$ with a smooth boundary, and that $f(z)$ has no zeros or poles on $\partial \Omega$. Then 
			$$ \dfrac{1}{2\pi i} \oint_{\partial \Omega} \dfrac{f'(z)}{f(z)} dz = N - P  $$ 
			Where $N$ and $P$ are integers that are equal to the number of zeros and poles of $f(z)$ in $\Omega$ respectively. 
		\end{theorem}
		The integer $|N-P|$ is also known as the {\em winding number} of the function $f(z)$. It is equal to the absolute value of the net number of times the image of $f(z)$ winds around the origin in $\C$ as the variable $z$ traverses the boundary $\partial \Omega$. 
		
		We apply this to the formula defining the Riccati-Evans functions in order to interpret the winding of the functions $E_{\bT}$ in terms of the roots of $D(\lambda)$. Suppose that we were in the chart corresponding to the matrix $\bT$. Denoting $\cdot := \dfrac{d}{d\lambda}$ we have 
		\begin{equation}\label{eq:ints}
		\begin{split}
		\oint_{\partial \Omega}\frac{ \dot{E}_{\bT}(\lambda)}{E_{\bT}(\lambda)} d \lambda & = \oint \frac{ \frac{d}{d\lambda} \left( \frac{D(\lambda)}{\det{\bX_\bT^u} \det{\bX_\bT^s}} \right)}{ \left( \frac{D(\lambda)}{\det{\bX_\bT^u} \det{\bX_\bT^s}} \right)} d \lambda  \\
		& = \oint \frac{\dot{D}(\lambda)}{D(\lambda)} d \lambda - \oint \frac{\det{\dot{\bX}_\bT^u}}{\det{\bX_\bT^u}} d \lambda 
		-\oint \frac{\det{\dot{\bX}_\bT^s}}{\det{\bX_\bT^s}} d \lambda
		\end{split}
		\end{equation}
		
		If we can choose a chart such that the $\det(\bX^{u,s}_\bT) \neq 0$ inside the simply connected domain $\Omega$, then the right two terms in \cref{eq:ints} vanish and the number of zeros of the Riccati-Evans function equals number of zeros of the original Evans function. 
		
		\section{(In)Stability Results: Application to the Model Equations} \label{sec:application}
		
		We apply the Riccati-Evans function described in \cref{sec:point} to first establish the numerical instability of travelling waves of type IV. We do this by tracking a real eigenvalue crossing zero into the right half plane as we lower the travelling wave speed below the minimal speed $c_*$ demarcating the transition from type II to type IV waves. We then numerically establish the stability of waves of type I, II and III by showing that for a reasonably large subset of the eigenvalue parameter $\lambda \in \C$, with $0 \leq \real{\lambda} \leq10^4$ there are no roots of the Evans function when $\bu$ is a travelling wave of speed $c>c_*$.
		
		We compute the Riccati-Evans function for \cref{eq:matslow} with asymptotic end states consisting of the stable subpace of $A_+$ and unstable subspace of $A_-$ for numerically computed waves of type I, II and IV. Without the precise wave speed of the type III waves, it is not possible to numerically solve for them, so all spectral data of the point spectrum must be inferred \cite{Harley_2014}. We used the continuation program AUTO to numerically compute travelling waves of type I, II and IV (and to approximate the minimal wave speed of type III), and used Mathematica's NDSolve function to solve the Riccati equation and compute the Riccati-Evans function. See \cref{fig:ric,,fig:type1real,fig:argbig,,fig:winding}.
		
		The only remaining ingredient is a (matrix for a) coordinate chart $\bT$. Finding such a chart can be a nontrivial task as there will inevitably be singularities in the matrix Riccati equation. The idea is to find a coordinate chart where the singularities do not appear in the region of the eigenvalue space we are interested in. For this system the matrix 
		$$
		\bT =  \begin{pmatrix} -i & 0 & 1 & 0 \\ 0 & i & 0 & 1 \\ 0 & 0 & i & 0 \\ 0 & 0 & 0 & -i \end{pmatrix} 
		$$
		was used and evidently produced no singularities of the Riccati equation (or the Riccati-Evans function) for values of $\lambda$ on the real line or in the upper right half of the complex plane (that we could observe numerically). A detailed determination of a chart that would always have this feature, as well as a proof of why that might be the case, is beyond the scope of this manuscript. 
		
		\subsection{Instability of type IV waves}
\begin{figure}[h]
			\makebox[1\linewidth][c]{%
				\begin{subfigure}[b]{0.47\textwidth}
					\centering
					\includegraphics[width=0.97\textwidth]{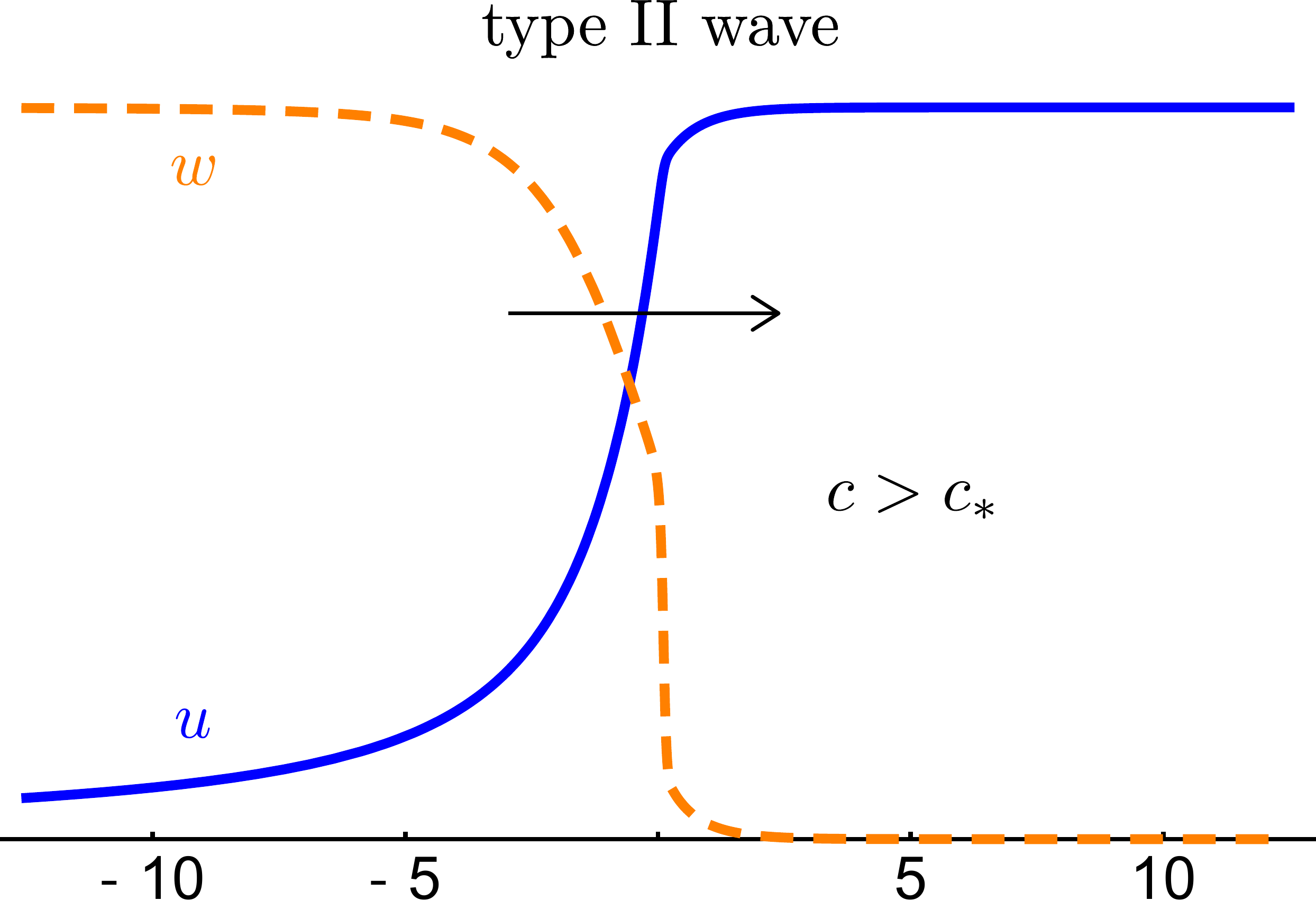}
				\end{subfigure}
				\begin{subfigure}[b]{.47\textwidth}
					\centering
					\includegraphics[width=0.97\textwidth]{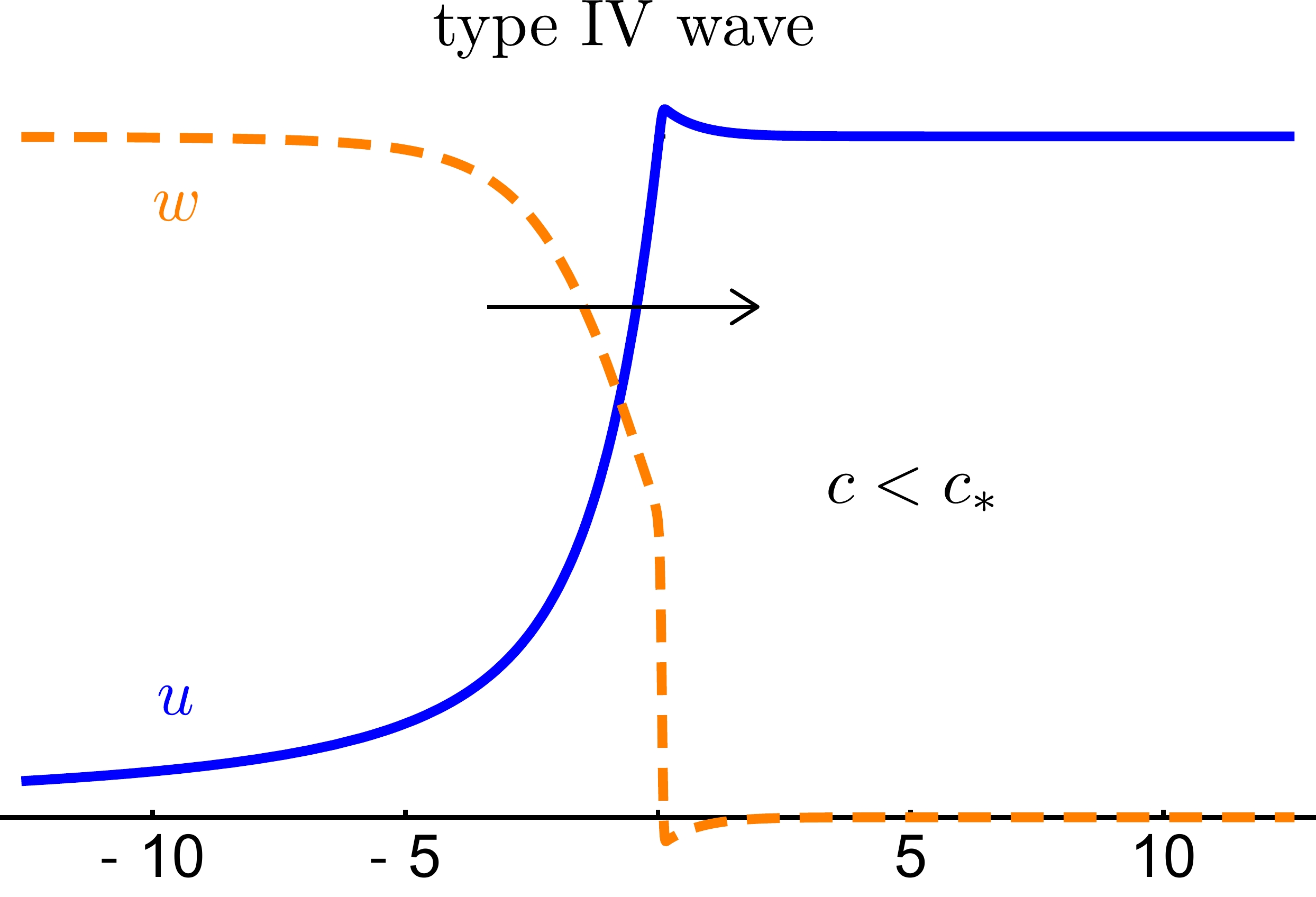}
				\end{subfigure}
			}
			\makebox[1\linewidth][c]{%
				\begin{subfigure}[b]{0.47\textwidth}
					\centering
					\includegraphics[width=0.97\textwidth]{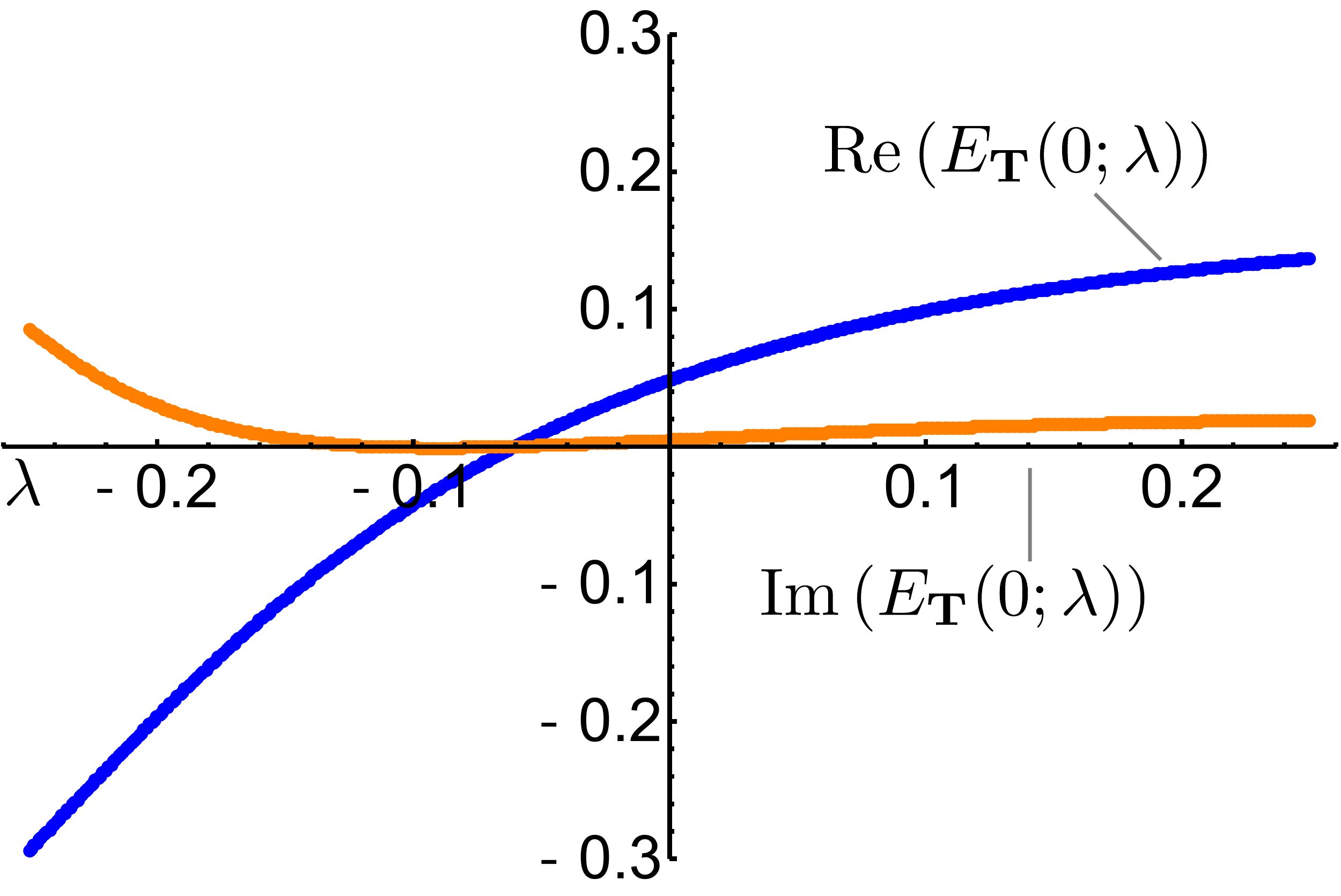}
				\end{subfigure}
				\begin{subfigure}[b]{.47\textwidth}
					\centering
					\includegraphics[width=0.97\textwidth]{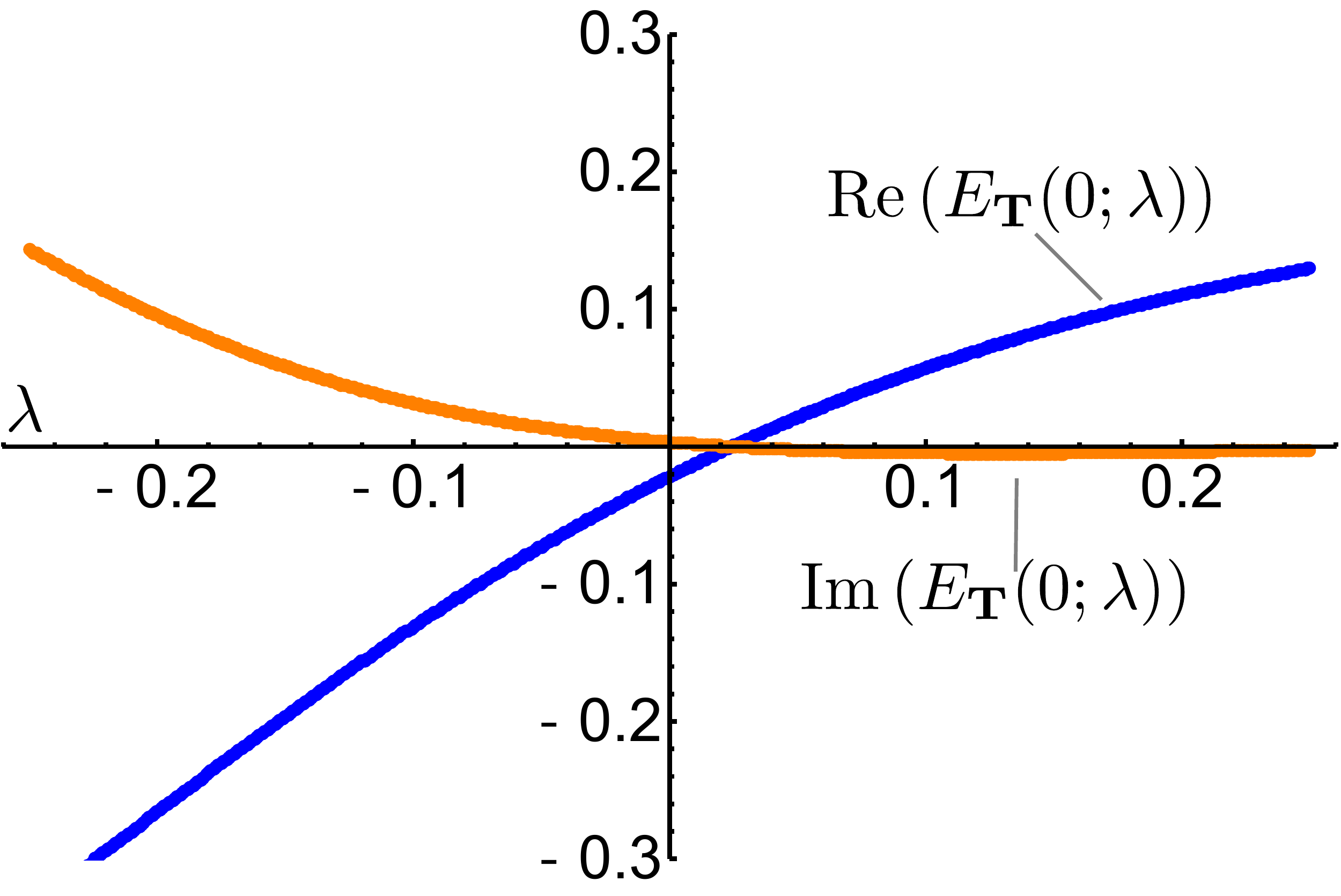}
				\end{subfigure}
			}
			\captionsetup{width=.95\linewidth}
			\caption{The top figures in each column show the type of wave that we are linearising about (Left column: type II, close to but slightly above the minimal wavespeed, and right column: type IV close to but slightly below). The bottom figures show the real and imaginary (blue and orange online respectively) parts of the Riccati-Evans function $E_\bT(0;\lambda)$, computed as a function of the (real) eigenvalue parameter $\lambda$. As the wave speed $c$ is decreased through the minimal wave speed ($c_*\approx 0.6701$), there is a real root of the Riccati-Evans function (and hence a real eigenvalue of the operator $L(\bu)$) which crosses into the right half plane, and as the type II waves transition to those of type IV, they become unstable.}
			\label{fig:ric}
		\end{figure}

		We first establish the instability of the type IV waves by plotting the Riccati-Evans function for real values of $\lambda$ and tracking a real eigenvalue as it crosses the imaginary axis as we lower the wave speed parameter $c$ below the threshold of the type III waves ($c_* \approx 0.6701$). See \cref{fig:ric}. From the plots of the Riccati-Evans function in the chart $\bT$, we see that for real values of $\lambda$ there do not appear to be any singularities of the function $E_\bT(0;\lambda)$, thus any zeros that appear are indeed zeros of the original Evans function and hence eigenvalues of the operator $L(\bu)$. There are many zeros on the real line, all of them negative until $c$ is made low enough, whereby the leading zero crosses into the right half plane. 
		
	\subsection{Stability of waves of type I, II and III}
\begin{figure}[ht]
			\centering
			\includegraphics[scale=0.37]{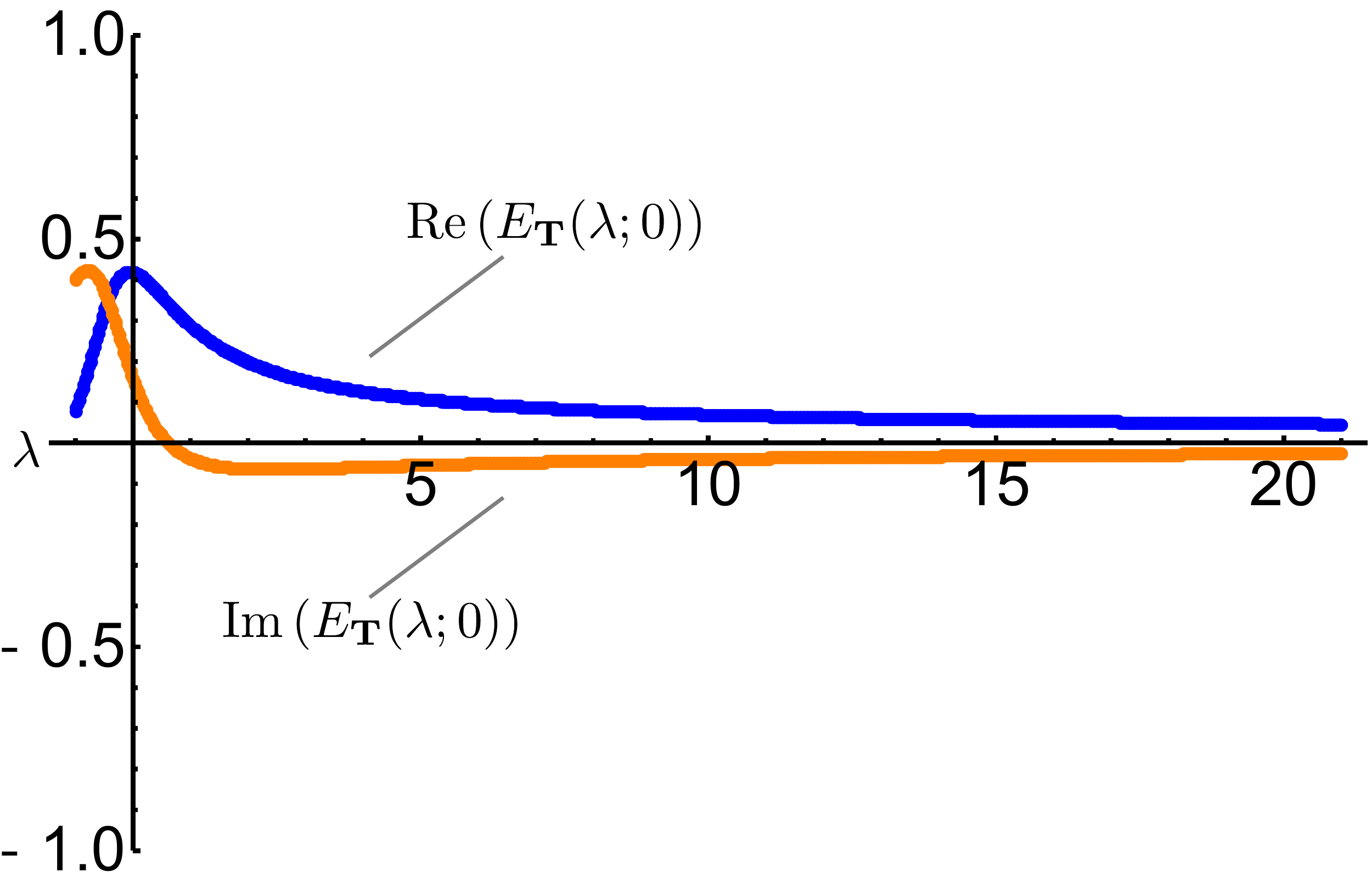}
			\captionsetup{width=.95\linewidth}
			\caption{A plot of the real and imaginary (blue and orange online) parts of the function $E_\bT$ for positive real values of the temporal spectral parameter $\lambda$ for the linearised operator about a type I wave. The parameter values are $u_\infty=1$, $c=1$ and $ \ve =0.01$.}
			\label{fig:type1real}
			\vspace{-.5cm}
		\end{figure}

\begin{figure}[ht]
			\centering
			\includegraphics[scale=0.29]{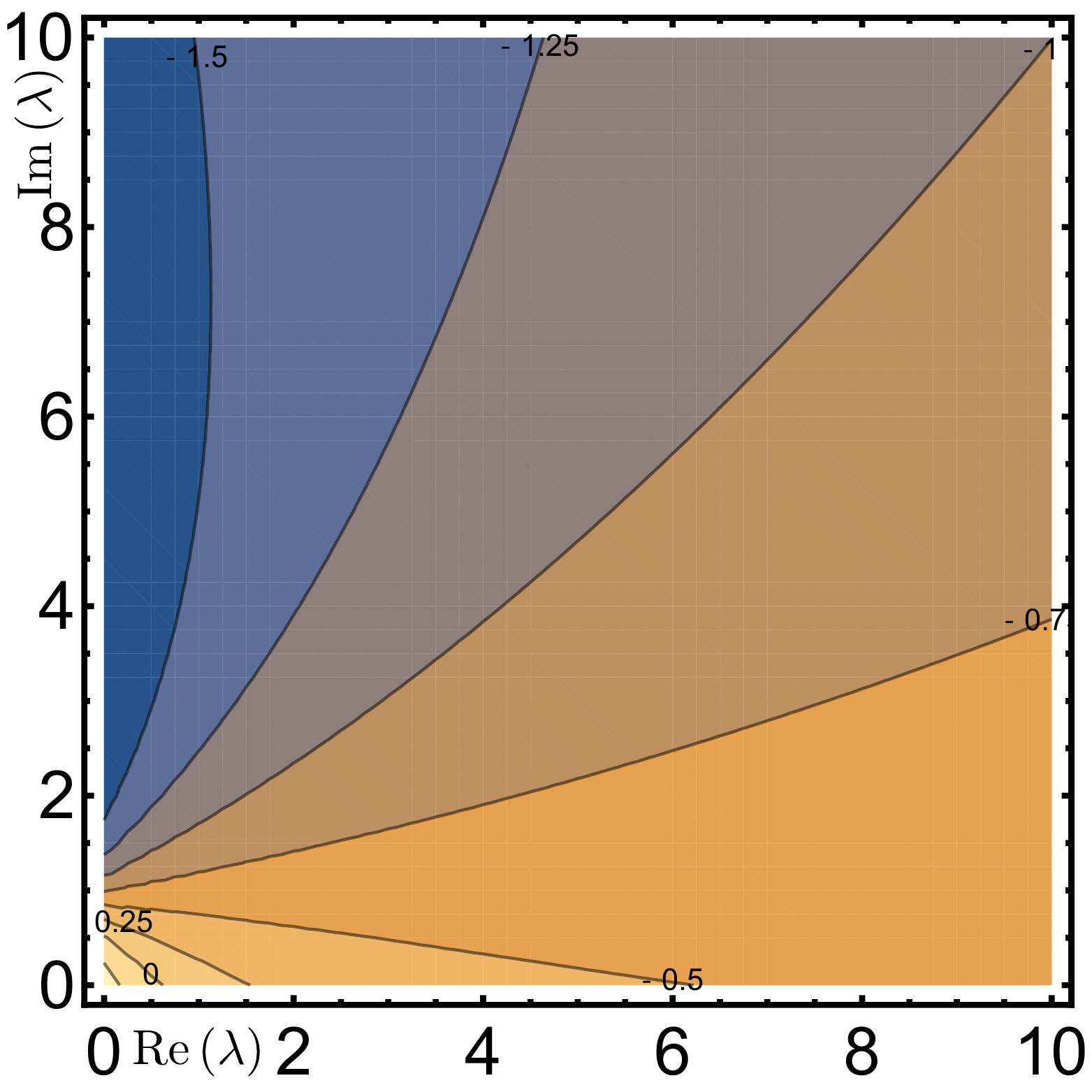} \quad
			\includegraphics[scale=0.35]{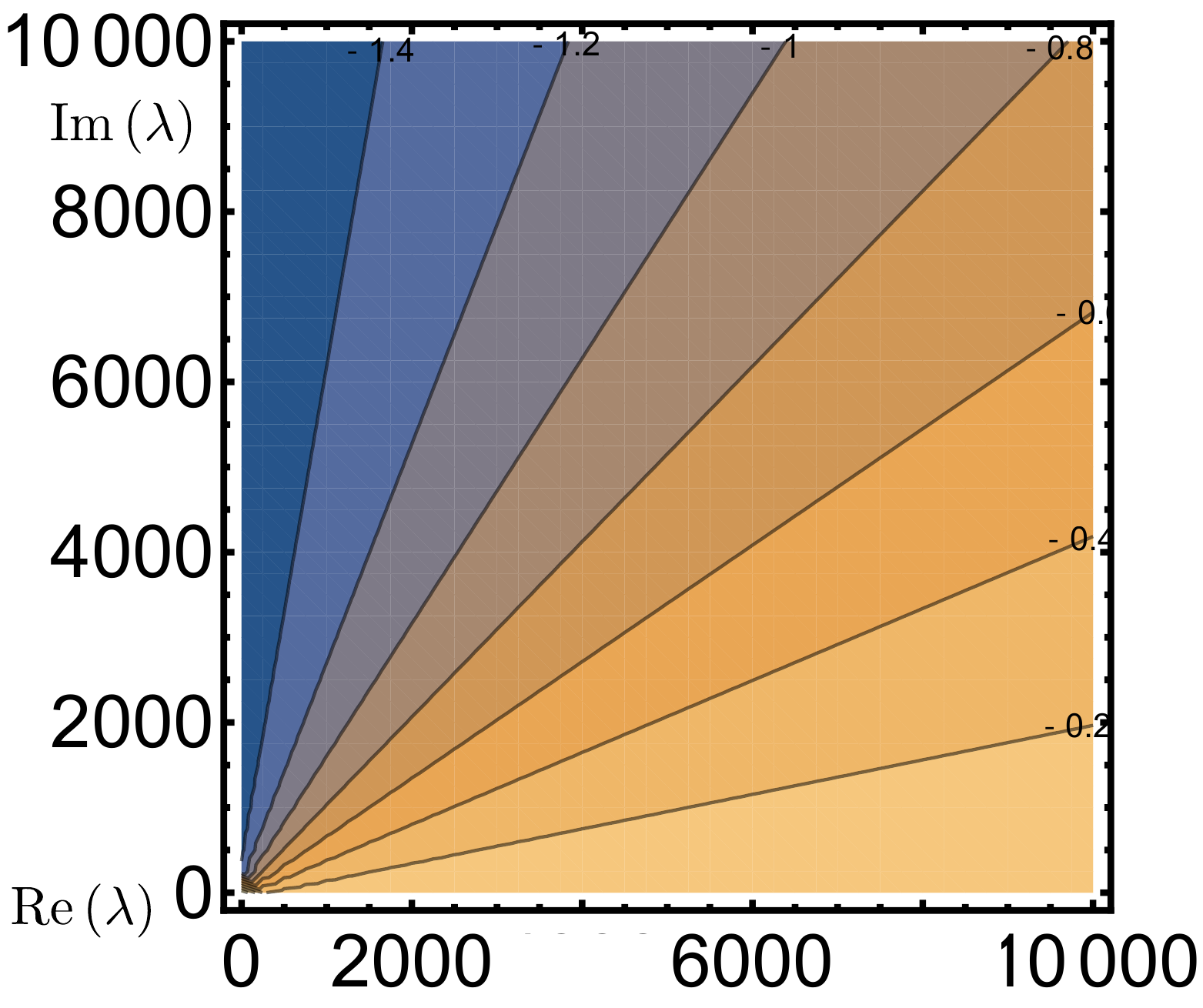}
			\caption{Left: A plot of contour lines of the argument of the function $E_\bT(\lambda;0)$ for the region of the first quadrant in the right half plane extending out to $\real{\lambda} <10$ and $\imag{\lambda} <10$. It is clear that there are no zeros or poles of the function $E_\bT$ in this region and hence no temporal eigenvalues. One can see the contour lines coalescing on a zero or a pole in the left half plane (in this case it is a pole). Right: A plot of contour lines of the argument of the function $E_\bT(\lambda;0)$ for the region of the first quadrant in the right half plane extending out to $\real{\lambda} <10,000$ and $\imag{\lambda} <10,000$. It is clear that there are no zeros or poles of the function $E_\bT$ in this region, and hence no temporal eigenvalues. Parameter values used were $u_\infty =1, c=1$, and $\ve = 0.01$. } 
			\label{fig:argbig}
			\vspace{-0.5cm}
		\end{figure}

\begin{figure}[h]
			\centering
			\includegraphics[scale=0.37]{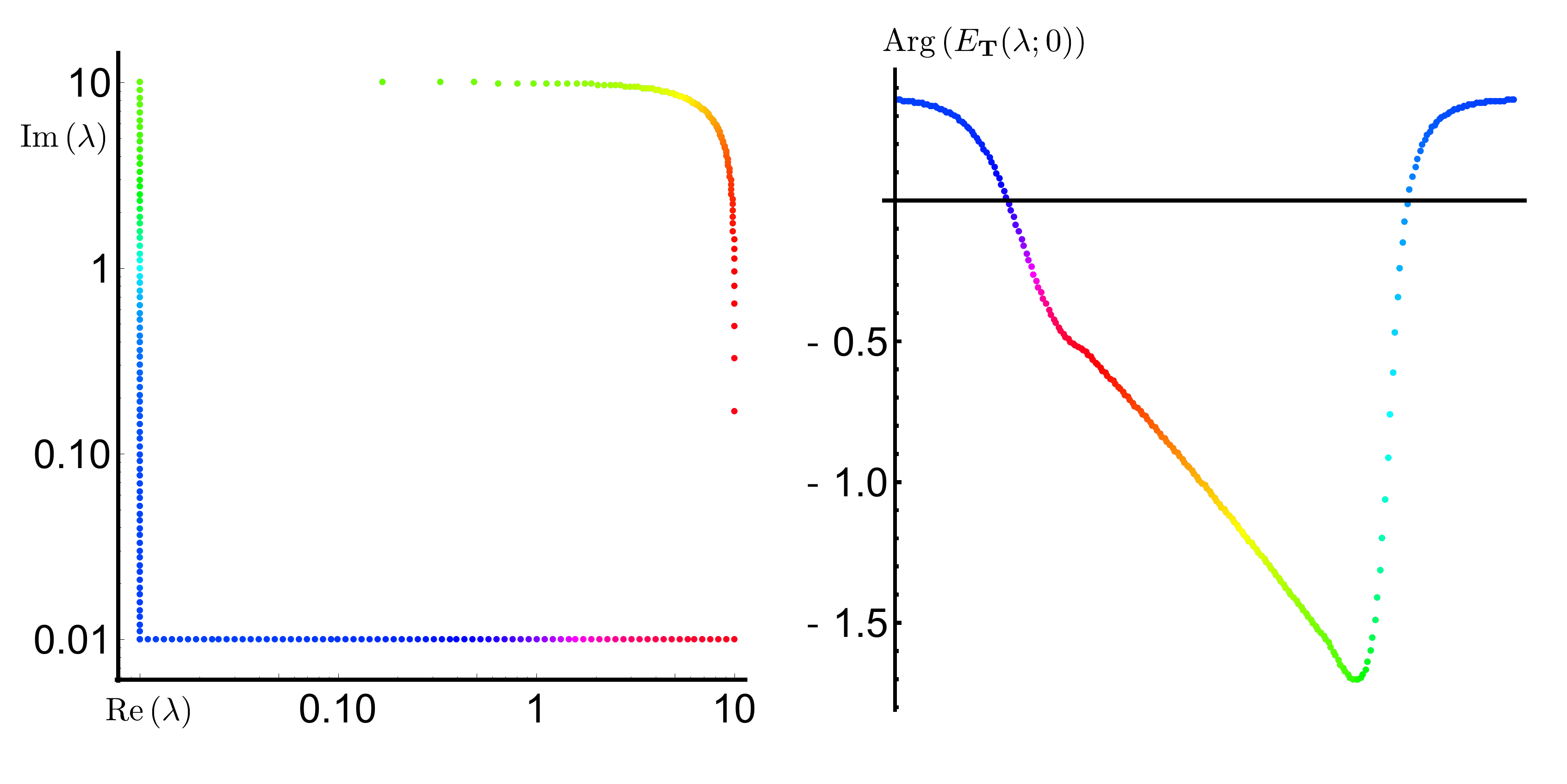} \,
			\includegraphics[scale=0.37]{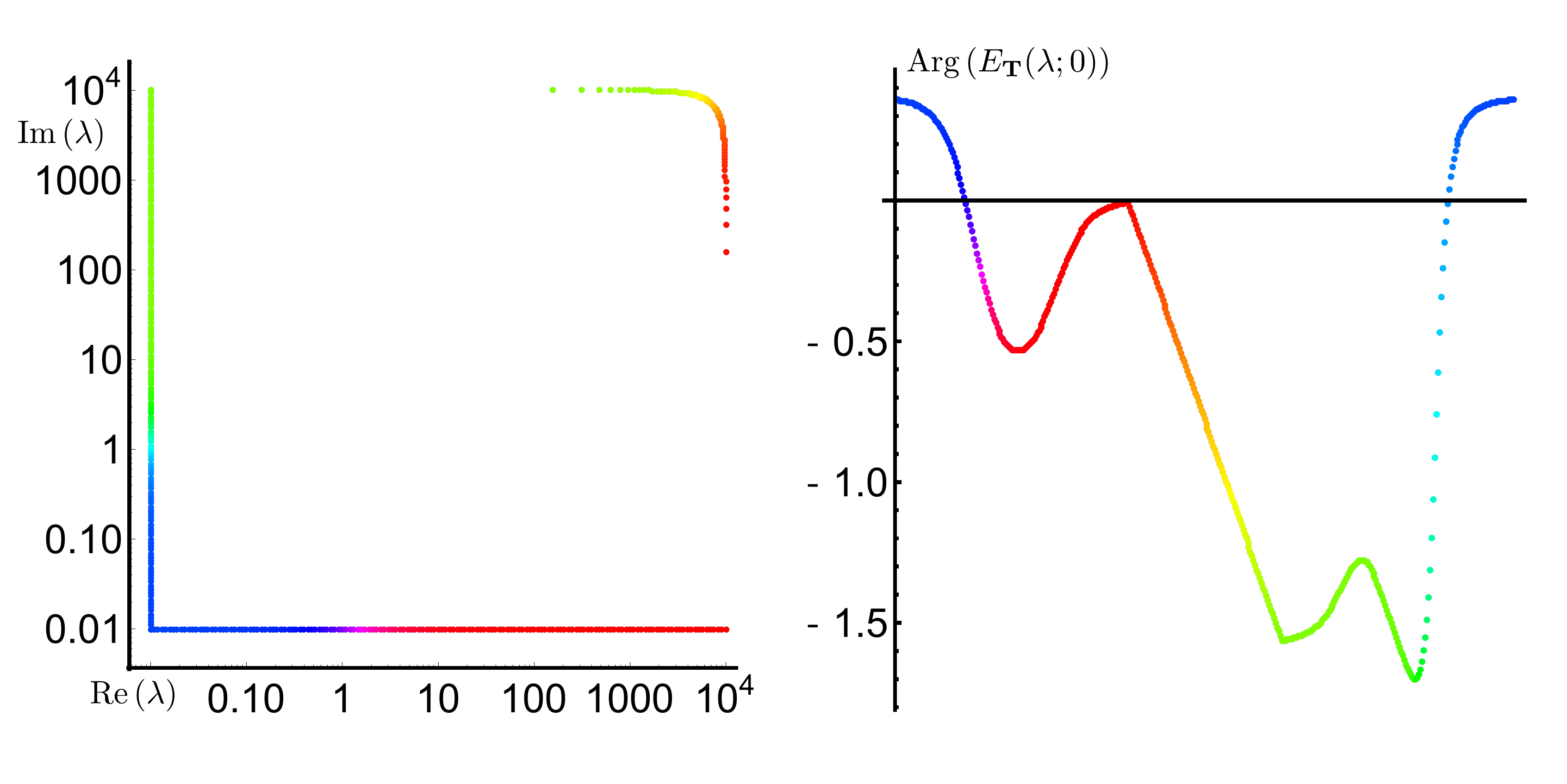}
			\captionsetup{width=.95\linewidth}
			\caption{(Colour online.) A plot of the function $\textrm{Arg}(E_\bT$) for values on the quarter circles of radius 10 (top) and 10,000 (bottom). The left figures depict a (logarithmic) parametrisation of the quarter circle, while the right figures are the corresponding plots (see colour online) of $\textrm{Arg}(E_\bT$). It is clear from the plots that there is no winding of the function $E_\bT$ here, and hence there is no spectrum of the linearised operator $L(\bu)$ in this region either. Parameter values used were $u_\infty =1, c=1$, and $\ve = 0.01$.} 
			\label{fig:winding}
			\vspace{-0.5cm}
		\end{figure}

		To numerically establish the spectral stability of travelling waves of type I and, II (and to infer spectral stability of the waves of type III), in the appropriately exponentially weighted spaces, we plot the argument of the Riccati-Evans function for successively larger regions in the upper right half plane. Because the travelling wave that we are linearising about is real, we know that any eigenvalues of the operator $L(\bu)$ must come in complex conjugate pairs, so if $\lambda$ is a root of $D(\lambda)$, then $\bar{\lambda}$ must also be a root of $D(\lambda)$. A consequence of \cref{eq:ric-evans-evans} is that, away from the poles of $E_\bT$, roots of the Riccati-Evans function must also come in conjugate pairs. Hence, it is sufficient to investigate the first quadrant of the complex plane for eigenvalues. In what follows, we show the numerical evidence for stability of type I waves only, the figures for waves of type II are qualitatively the same. \Cref{fig:type1real} shows a plot of the function $E_\bT(\lambda;0)$ for real values of $\lambda$. It is clear that there are no roots of the Riccati-Evans function for $\lambda <20$. To investigate complex eigenvalues, we plot the argument of the function $E_\bT$ a large section of the complex plane. For a meromorphic function, a zero or a pole is represented by the coalescing of many contour lines of the argument of the function. Hence, we can visually see from \cref{fig:argbig} that there are no zeros or poles of the linearised operator $L(\bu)$ for the type I wave in this region of $\C$. We confirm this with the argument principle by computing the winding number of the Riccati-Evans function on successively larger quarter circles and can again visually see that no winding takes place (see \cref{fig:winding}). 

		\section[Future work]{Discussion and future work}\label{sec:fut} 
In this manuscript, we studied the spectral stability of the four different types of travelling waves supported by an advection-reaction-diffusion equation originally proposed in \cite{PER} to describe haptotactic cell invasion in a model for melanoma. Using a Riccati-Evans function approach, we numerically showed that the biologically-unfeasible type IV waves -- waves for which the invasive tumour cell population wave profile $w$ is negative for certain parts of the profile -- are unstable, while the other three types of waves where the tumour cell population $w$ stays positive are spectrally stable. Heuristically, instability of the type IV waves follows from the fact that the type III waves have a (very) fast decay at $+\infty$. Thus $\lambda =0$, the eigenvalue associated with spatial invariance of the front, is a temporal eigenvalue in the now weighted space. It persists, and in this case moves into the right half-plane as the wave-speed is further decreased (which is what we numerically showed).

A logical next step is to further study the connection between the observed wave speed and the asymptotic behaviour of its initial condition.
This connection was already partly investigated in \cite{Harley_2014, perumpanani2000traveling}. In \cite{Harley_2014}, formal computations around the asymptotic end state of a travelling wave are used to show that the type I and type II waves travel with speed $c= 1/\chi + \mathcal{O}(\ve)$, where $\chi$ is the asymptotic decay rate at $\infty$ of the exponentially decaying initial condition for $w$ (i.e. $w(x,0) =w_0(x) = \max\{1,e^{-\chi x}\}$). This result was also numerically verified in \cite{perumpanani2000traveling}. Unfortunately, the asymptotic linear analysis of \cite{Harley_2014} was unable to derive a correct approximation for the minimal wave speed $c_*$ associated with the type III waves (i.e. the type III waves are pushed fronts \cite{van2003front}),
see in particular  \cite[Fig. 10]{Harley_2014}. In \cite{perumpanani2000traveling}, the authors used a power series approximation to derive a quadratic relationship between the minimal wave speed $c_*$ and the asymptotic end state of the wave $u_{\infty}$  in the singular limit $\ve=0$. Combining the results of \cite{Harley_2014} and \cite{perumpanani2000traveling} indicated that $c_* = c_*(u_\infty,\ve)$ and it remains to be seen if this relationship can be derived analytically.

We are currently working on using this approach to study the stability of travelling waves in a model for wound healing angiogenesis \cite{H1}, a model for {stellar wind} \cite{carter2017transonic}, and in two different types of tumour invasion models \cite{H3, H2}. The Riccati-Evans function approach in this manuscript does not take advantage of the singularly perturbed nature of the stability problem. The nonlocal eigenvalue problem approach \cite{DGK2,DGK1,HDK08} and the singular limit eigenvalue problem approach \cite{nishiura1987stability, nishiura1990singular} are two related analytical techniques that use this singular perturbed nature to simplify the Evans function computations. It would be interesting to see if, similar to \cite{de2016spectra}, one of these techniques can be incorporated in the Riccati-Evans function approach to further optimise the computations. In particular, in \cite{de2016spectra} the authors use the Riccati equation and the singularly perturbed nature of the problem to compute a factored Evans function via the Grassmanian, where one of the factors is analytic and never zero, thus reducing the calculations necessary for eigenvalue determination. We comment that the factorisation of the Evans function given by \cref{eq:ric-evans-evans} is reminiscent of that in \cite{de2016spectra} (when the chart is chosen properly) - though it does not make use of any singular structure in the problem.
		
We note that in \cite{grudzien2016geometric}, the authors factor the Evans function in a different way, reducing the computations to ones in a unitary matrix (Hopf) bundle. The factorisation in \cref{eq:ric-evans-evans} is seemingly complementary to that in \cite{grudzien2016geometric} in the sense that the unstable bundle in \cite{agj90} factors into two sub-bundles, the transition maps of one being the unitary group, while the transition maps of the other are the Grassmannian (in the sense that it is a homogeneous space of Lie groups). 
		  
\section*{Acknowledgements} The authors would like to thank G. Gottwald, D. Lloyd, and A. G. Munoz for their helpful numerical advice, as well as the referees for their valuable input and suggestions. RM would like to thank S. J. Malham and M. Beck for very insightful conversations regarding the Grassmannian of two planes in $\C^4$ and RM and TVR would like to thank D. Smith for his commentary on the argument principle in complex analysis. PvH acknowledges support under the Australian Research Council grant DE140100741. MW acknowledges support under the Australian Research Council grant DP180103022.

		\bibliographystyle{amsalpha}
		\bibliography{references3}

\end{document}